\documentclass[12pt]{amsart}

\usepackage{amssymb}

\usepackage{enumitem}

\usepackage{graphicx}

\makeatletter
\@namedef{subjclassname@2020}{%
  \textup{2020} Mathematics Subject Classification}
\makeatother

\usepackage[T1]{fontenc}

\newcommand{\BN}{{\mathbb N}}
\newcommand{\BR}{{\mathbb R}}

\newcommand{\FF}{{\mathcal{F}}}

\newcommand{\BB}{{\mathcal{B}}}

\newcommand{\MM}{{\mathcal{M}}}

\newcommand{\EE}{{\mathcal{E}}}

\newtheorem{theorem}{Theorem}[section]
\newtheorem{corollary}[theorem]{Corollary}
\newtheorem{lemma}[theorem]{Lemma}
\newtheorem{proposition}[theorem]{Proposition}

\theoremstyle{definition}
\newtheorem{definition}[theorem]{Definition}
\newtheorem{remark}[theorem]{Remark}
\newtheorem{example}[theorem]{Example}

\numberwithin{equation}{section}

\begin{document}


\baselineskip=17pt


\title[Dirichlet forms and the obstacle problem]{Quasi-regular Dirichlet forms and the obstacle problem for
elliptic equations with  measure data}

\author[T. Klimsiak]{Tomasz Klimsiak}
\address{Faculty of Mathematics and Computer Science\\ Nicolaus Copernicus University\\
Chopina 12/18\\
87-100 Toru\'n, Poland}
\email{tomas@mat.umk.pl}

\date{}

\begin{abstract}
We consider the obstacle problem with irregular barriers for
semilinear elliptic equations involving measure data and operator
corresponding to a general quasi-regular Dirichlet form. We prove
existence and uniqueness of a solution as well as its
representation as an envelope of  a supersolution to some related
partial differential equation. We also prove regularity results
for the solution and the Lewy-Stampacchia inequality.
\end{abstract}



\maketitle

\section{Introduction}
Let $E$ be a Lusin space (i.e. the image of a Polish space under injective continuous mapping), $m$ be a $\sigma$-finite positive
measure on $\mathcal{B}(E)$ and let $(L,D(L))$ be a Dirichlet
operator associated with some quasi-regular (possibly non-symmetric)
Dirichlet form $(\EE,D[\EE])$ on $L^2(E;m)$. In the present paper,
we investigate the obstacle problem of the form
\begin{equation}
\label{eq1.1} \left\{
\begin{array}{l}-Lu\le f(\cdot,u)+\mu \quad
\mbox{on }\{u>h_{1}\},\smallskip \\
-Lu\ge f(\cdot,u)+\mu \quad\mbox{on }\{u<h_{2}\}, \smallskip \\
h_{1}\le u\le h_{2} \quad m\mbox{-a.e.},
\end{array}
\right.
\end{equation}
where $\mu$ is a smooth measure (if $\mu$ is bounded this means that $\mu$ charges no $\EE$-exceptional sets; for general definition see Section \ref{sec2}), $f:E\times\BR\rightarrow\BR$ and $h_{1},h_{2}$
are measurable functions on $E$ such that $h_{1}\le h_{2}$
$m$-a.e. We also consider one-sided problem, i.e. we allow
$h_1\equiv-\infty$ or $h_2\equiv+\infty$.

The class of operators associated with quasi-regular Dirichlet
forms is quite wide. It includes local operators in divergence
form,  $\alpha$-Laplacian type operators,
Ornstein-Uhlenbeck type operators in Hilbert spaces and others
(see, e.g., \cite{Fukushima,J1,KR:JFA,KR:CM,MR} for concrete
examples). We think that the fact that nonlocal operators fit into
our general framework is of special interest, because problem
(\ref{eq1.1}) with nonlocal operators and measure data is
considered here for the first time.

For an overview of numerous applications of
obstacle problem for elliptic and parabolic PDEs  we  refer the reader to \cite{XRO} and references
therein. In recent years nonlocal models attracted quite a lot of
interest because  it was observed that equations with nonlocal
L\'evy type operators describe some phenomena better then local
equations (see, e.g., \cite{CT,HQD}). The applications we have in
mind include population biology models, models of mathematical
finance involving jump processes and some interacting particles
models with repulsive/attractive interaction potentials. In all
the mentioned models the obstacle problem with rough data
naturally appears. In population models it is well known (see
\cite{DD1,DHMP})  that solutions of steady-state predator-pray
models with sufficiently large/small appropriate parameters behave
like solutions of certain free boundary problems which may be
equivalently formulated as an obstacle problem of the form
(\ref{eq1.1}) with merely measurable barrier. In these models, $L$
describes the dispersal of animals, $f$ describes the growth of
population and $\mu$ is the harvesting distribution. In the theory
of option pricing, the fair price of some derivative contracts are
of the form
\begin{align}
\label{eq1.00.1} u(x)&=\sup_{\tau\ge 0}\inf_{\sigma\ge 0}
E_x\Big(\int_0^{\tau\wedge\sigma}f(\cdot,u)(X_r)\,dr\nonumber\\
&\quad+\int_0^{\tau\wedge
\sigma}dA^\mu_r+h_1(X_\tau)\mathbf{1}_{\{\tau<\sigma\}}
+h_2(X_\sigma)\mathbf{1}_{\{\sigma\le
\tau\}}\Big),
\end{align}
where $X$ is a process with generator $L$ starting from $x$ at
time 0. The process $X$ describes  the evolution of stock prices,
$f$ generates the nonlinear expectation (see \cite{DQS}), the
additive functional $A^\mu$ (generated by a smooth measure $\mu$)
is the running  cost or profit, and $h_1(X), h_2(X)$ are pay-off
processes (such a situation appears for instance when considering
American options or Israeli options). Since 70', connections of
value functions of the form (\ref{eq1.00.1}) with obstacle
problems with one and two barriers have been intensively studied
in the literature (see, e.g., \cite{BL,GIOQ,Zabczyk}). It is worth
noting here that in some applications (for instance to digital
options, see \cite{GIOQ}) the functions $h_1, h_2$ are assumed to
be merely measurable. In the interacting particles models,  for
given Green function  $G$ and positive function $W_a$, we are
looking for a local minimizer for the interacting energy
\begin{equation}
\label{eq1.00.2}
E[\gamma]=\frac12\int_E\int_E(G(x,y)+W_a(x-y))\gamma(dy)\,\gamma(dx)
\end{equation}
in the class of  probability measures $\gamma$ on $E$. It is known
(see, e.g., \cite{CDM} for the case of Riesz's potentials) that
such a minimizer is the second component (see (\ref{eq1.2})) of
the local solution to (\ref{eq1.1}) with operator $L$ associated
with the Green function $G$ and $\mu=-LW_a* \gamma$.

In this paper, we impose very weak conditions on $\EE$ and  the data
$\mu$, $f$, $h_1$, $h_2$. To formulate them,  let us first recall that the operator
$(L,D(L))$
and its adjoint operator $(\hat L,D(\hat L))$ are generators of $C_0$-semigroups
of contractions $\{T_t,\, t\ge 0\},\, \{\hat T_t,\, t\ge 0\}$ on $L^p(E;m)$ for every $p\ge 1$.
Let $\{G_\alpha,\, \alpha>0\}$ (resp. $\{\hat G_\alpha,\, \alpha>0\}$) be the resolvent of $\{T_t,\, t\ge 0\}$
(resp. $\{\hat T_t,\, t\ge 0\}$). For positive $f\in L^p(E;m)$ we set
\[
Gf=\sup_{n\ge 1}G_{1/n}f,\quad \hat Gf=\sup_{n\ge 1}\hat G_{1/n}f.
\]
In the paper, we assume that $\EE$ satisfies strong sector condition (see Section 2) and it is transient, i.e. $Gf$ is finite $m$-a.e. for some strictly
positive $f\in L^1(E;m)$ (and hence for every $f\in L^1(E;m)$). It is known (see \cite{Fukushima}) that this condition is
 equivalent to the existence of  a strictly positive $g\in L^1(E;m)$ such that
\begin{equation}
\label{eq1.10}
\int_E|u|g\,dm\le\EE(u,u)^{1/2},\qquad u\in D[\EE].
\end{equation}
As for $\mu$, we assume that it belongs to the class
\begin{align}
\label{eq1.02} \mathbb M_0&=\{\mbox{$\mu:|\mu|$ is smooth and
$\hat{G}\phi\cdot\mu\in\mathcal{M}_{0,b}$ for}\nonumber\\
&\qquad\qquad \mbox{some $\phi\in L^1(E;m)$ such that $\phi>0$
$m$-a.e.}\}
\end{align}
considered in \cite{KR:CM}. Here $|\mu|$ denotes the variation of
$\mu$, $\mathcal{M}_{0,b}$ is the space of all finite smooth
signed measures on $\BB(E)$. Of course, the class $\mathbb M_0$
depends on the structure of $\EE$, but by \cite[Corollary
1.3.6]{O}, we always have $\MM_{0,b}\subset\mathbb M_0$. In
general, the inclusion is strict. For instance, if $d\ge3$ and
$L=\Delta^{\alpha/2}$ with $\alpha\in(0,2]$ on a bounded smooth
domain $D\subset\mathbb{R}^{d}$, then by \cite{Ku} there exist
$c_1,c_2>0$ such that
\begin{equation}
\label{eq1.11}
c_2\delta^{\alpha/2}(x)\le G1(x)\le c_2\delta^{\alpha/2}(x),
\quad x\in D,
\end{equation}
where $\delta$ is the distance to the boundary of $D$.
Consequently, in that case $\mathbb M_0$ includes Radon measures
of infinite total variation. In particular, we have
$L^{1}(D;\delta^{\alpha/2}(x)\,dx)\subset \mathbb M_0$. In recent
years elliptic equations involving Laplace operator and $L^{1}(D;
\delta(x)\,dx)$ data were considered by many authors (see, e.g.,
\cite{QS,Rakotoson} and references therein). Note that it also may
happen that $\mathbb M_0$ includes nowhere Radon measures (see
Example \ref{ex.ex}). If the resolvent of  the operator $(L,D(L))$
is strongly Feller (i.e. $G_\alpha (\BB_b(E))\subset C_b(E)$), then
$(L,D(L))$ has the Green function, i.e. there exists $r\in
\BB^+(E)\times\BB^+(E)$ such that
\[
Gf=\int_Er(\cdot,y)f(y)\,m(dy),\quad f\in L^1(E;m),
\]
and moreover,
\[
\mathbb M_0\supset\{\mu\mbox{ is a Borel measure on }E:
\int_Er(x,y)|\mu|(dy)<\infty,\,x\in E\}.
\]
The inclusion above  can be replaced by equality if we additionally assume  that
$\mu$ is smooth and replace ``for every" by ``quasi every" (with
respect to the capacity associated with $\EE$). The
characterization of $\mathbb M_0$ in this spirit is also possible
for general operator $(L,D(L))$ but to state it requires the
introduction of  the notion of positive additive functional (see
Section 3).

The function $f: E\times\mathbb{R}\rightarrow\mathbb{R}$ is
assumed to be continuous and nonincreasing with respect to the
second variable. We also assume that $f(\cdot,0)\cdot m \in\mathbb
M_0$ and for every $y\in\BR$, $f(\cdot,y)$ is quasi-integrable
(weaker condition than integrability, see Section \ref{sec2}). These assumptions on $f$ were used previously in
many papers devoted to linear and nonlinear equations involving
measure data and local operators but with $f(\cdot,y)\in L^1(E;m)$
(see, e.g., \cite{Betal.,Boccardo}). Semilinear elliptic equations
with quasi-integrable data and local operators were considered in
\cite{OP}. Equations with quasi-integrable data and nonlocal
operators were considered for the first time in \cite{KR:JFA} (see
also \cite{KR:CM}).

In the paper we do not impose any regularity assumption on the
barriers $h_1,h_2$.  Therefore to guarantee the existence of a
solution we have to assume that they satisfy some kind of
separation condition. Roughly speaking, our  condition says  (see
Section 3) that between the barriers one can find some function
$v$ such that $v$ is a difference of two natural potentials and
$f(\cdot,v)\cdot m\in \mathbb M_0$. For instance,
this condition is satisfied if $h_1\le \varphi(w)\le h_2$ and
$f(\cdot,\varphi(w))\in L^1(E;m)$ for some $w\in D(L)$ and
$\varphi$ being  difference of two convex function and such that
$\varphi(0)=0$.

Since our data are irregular, the classical approach to
(\ref{eq1.1}) via variational inequalities (see \cite{AP,Brezis,
Stampacchia}) does not apply (see, however, \cite{BrezisSerfaty}
for the case $L=\Delta$). In the present paper by a solution to
(\ref{eq1.1}) we understand  a solution of the complementary
system (see \cite{AP,KS}) associated with (\ref{eq1.1}). Roughly
speaking it is a pair $(u,\nu)$ consisting of a quasi-continuous
function $u$ on $E$ and a measure $\nu\in\mathbb M_0$ such that
\begin{equation}
\label{eq1.2} \left\{
\begin{array}{l}-Lu= f(\cdot,u)+\mu+\nu, \smallskip \\
h_{1}\le u\le h_{2} \quad m\mbox{-a.e.}, \smallskip \\
\int_E(u-h_{1})\,d\nu^{+}=\int_E(h_{2}-u)\,d\nu^{-}=0,
\end{array}
\right.
\end{equation}
where $\nu^+,\nu^-$ denote the positive and negative parts in the
Jordan decomposition of $\nu$.

The obstacle problem with irregular data is a subject of intensive
study. Most of available results are formulated in the language of
differential inclusions (when $L$ is a  general accretive or
completely accretive operator) or in the language of entropy or renormalized
solutions (when $L$ is a nonlinear Leray-Lions type
operator; when $L$ is a linear Leray-Lions type operator, one can use an equivalent notion of Stampacchia's solution by duality).

The paper by Brezis and Strauss \cite{BS} is the first paper
devoted to problem of type (\ref{eq1.1}) with $L^1$ data. More
precisely, in \cite{BS} differential inclusions of the form
\begin{equation}
\label{eq1.inc} -\lambda u-Au+\beta(x,u)\ni\mu
\end{equation}
are considered. In (\ref{eq1.inc}), $\lambda\ge0$, $\mu\in
L^1(E;m)$, $A$ is an operator with sub-Markovian resolvent such that $D(A)\subset L^1(E;m)$,
and for fixed $x\in\BR$, $\beta(x,\cdot)$  is a maximal monotone
graph on $\BR\times\BR$. Note that if we define $\beta$ by
\begin{equation}
\label{eq1.mmg} D(\beta(x,\cdot))=[h_1(x),h_2(x)],\quad
\beta(x,y)=\left\{
\begin{array}{ll} [0,\infty),\quad& y=h_1(x),\\
\{0\},\quad & h_1(x)<y<h_2(x),\\
(-\infty,0],\quad & y=h_2(x),
\end{array}
\right.
\end{equation}
then (\ref{eq1.inc}) reduces to the obstacle problem with operator
$L=\lambda+A$ and barriers $h_1$ and $h_2$. In fact, in \cite{BS}
equation (\ref{eq1.inc}) with $\lambda=0$ and $\beta$ not
depending on  $x$ is considered, so the results of \cite{BS}
apply to obstacle problems with constant barriers.  As for $A$, in
\cite{BS} it is assumed that
\begin{equation}
\label{eq1.l1} \|u\|_{L_1}\le c\|Au\|_{L^1},\quad u\in D(A).
\end{equation}
The above conditions guarantee that the solution $u$ to
(\ref{eq1.inc}) belongs to $D(A)\subset L^1(E;m)$.
Consequently, if we set
\begin{equation}
\label{eq1.sub} w:=\mu+\lambda u+Au,
\end{equation}
then $w\in L^1(E;m)$ (of course $w\in\beta(u)$ a.e.). By the
monotonicity of $\beta$, for every function $v$ on $E$ such that
$h_1\le v\le h_2$, we have
\[
\int_E(u-v)w\,dm =\int_E(u-v)(w-0)\,dm\le 0
\]
since $0\in\beta(v)$ a.e. In different words, the pair $(u,w\cdot m)$ is a solution to (\ref{eq1.2}) with $L=\lambda+A$.

When $\beta$ depends on $x$, then depending on the regularity of
$\beta$ with respect to $x$, one can consider  the so called strong or
generalized solutions to (\ref{eq1.inc}). Hence, in the case where
$\beta$ is given by (\ref{eq1.mmg}), the concept of  solution
depends on the regularity of barriers (see \cite{Wittbold,Voigt}).
Roughly speaking, strong solution corresponds to the case, when
the reaction measure $\nu$ (or, equivalently, $w$) is absolutely
continuous with respect to $m$. Generalized solutions to
(\ref{eq1.inc}) with $\mu\in L^1(E;m)$ were considered in
\cite{BW,Wittbold}.  In \cite{Wittbold}
problem (\ref{eq1.inc}) with  a  linear  Leray-Lions type operator $A$ is considered. It is shown there that in general $w$ is a measure and for every function $v$  on $E$ such that
$h_1\le v\le h_2$,
\[
\int_E(u-v)\,dw\le 0.
\]
Therefore also in case $\beta$ depends on $x$ problem (\ref{eq1.inc}) can be rewritten in the form (\ref{eq1.2}) (see also \cite[Theorem 3.2]{AP}).

The obstacle problem of the form (\ref{eq1.1})  with a   nonlinear
Leray-Lions type operator $L$ and $\mu\in L^1(E;m)$ was considered
in \cite{BoccardoCirmi,BoccardoGallouet}. In both papers the
problem is studied in the setting of  entropy solutions introduced
in \cite{Betal.}  (for a closely related notion of renormalized
solution see \cite{DMOP}).

To our knowledge, first results concerning  (\ref{eq1.1}) with
``true" measure data were  obtained in \cite{DL} by using
Stampacchia's approach by duality (see also \cite{DD}). In
\cite{DL} the obstacle problem with one lower barrier $h_1$ (i.e.
$h_2\equiv +\infty$) is considered and it is assumed that $L$ is a
uniformly elliptic divergence form operator. The results of
\cite{DL} were  extended  in \cite{Leone} to the case of nonlinear
Leray-Lions type operator $L$.  In \cite{Leone}  the setting of
renormalized solutions is used.

Quite recently first papers devoted to semilinear elliptic
equations involving measure data and nonlocal operators  (mostly
fractional Laplacian) appeared (see, e.g., \cite{AAB,CV,KPU,KMS}).
General results on existence, uniqueness and regularity of
solutions of such equations with operator $L$ corresponding to
Dirichlet form were proved in \cite{KR:JFA,KR:CM} (see also
\cite{KR:NODEA}) in case $\mu$ is a smooth measure, and in
\cite{Kl:CVPDE} for a general Borel measure $\mu$. However,  to
our knowledge, there are no results on obstacle problem
(\ref{eq1.1}) with true measure data and nonlocal operator $L$.
Therefore all the results of the present paper are new in case $L$
is nonlocal and $\mu$ is a ``true" measure. It is worth
mentioning, however, that they are new even if $\mu\in L^1(E;m)$,
because as  compared with papers devoted to  problem
(\ref{eq1.inc}) we consider the case $\lambda=0$ and we do not
assume (\ref{eq1.l1}). Also note that in general, our solutions
are not even locally integrable, so need not satisfy the condition
\[
\int_E(u-k)^+\,dm<\infty\quad\mbox{for some}\quad k>0,
\]
which is the minimal requirement on $u$  when one investigates (\ref{eq1.inc}) in the setting of completely accretive operators
(see \cite{BC}).

In general, under  weak assumptions on $f,\mu$ described above
the solution $u$ to  (\ref{eq1.1}) may be very irregular.
Therefore the problem of making sense of the first equation in
(\ref{eq1.2}) arises. Following \cite{KR:JFA,KR:CM} we address it
by using stochastic analysis methods. Namely, by a solution of the
first equation in (\ref{eq1.2}) we mean a function
$u:E\rightarrow\BR$ satisfying for quasi-every (q.e. for short)
$x\in E$ the following generalized Feynman-Kac formula
\begin{equation}
\label{eq1.3} u(x)=E_{x}\int_{0}^{\zeta}f(X_{t}, u(X_{t}))\,dt
+E_{x}\int_{0}^{\zeta}dA_{t}^{\mu}+
E_{x}\int_{0}^{\zeta}dA_{t}^{\nu}.
\end{equation}
Here $\mathbf{M}=(X,P_{x})$ is a  special standard process with
life-time $\zeta$ associated with the form $(\EE,D[\EE])$, $E_x$
is the expectation with respect to $P_x$ and $A^{\mu}, A^{\nu}$
are continuous additive functionals of $\mathbf{M}$ in the Revuz
correspondence  with $\mu$ and $\nu$, respectively.

It is worth remarking that in the important case where $\mu,\nu\in
\mathcal{M}_{0,b}$, the probabilistic definition (\ref{eq1.3}) can be
rephrased in purely analytical terms. Namely, under these
assumptions on $\mu,\nu$, (\ref{eq1.3}) is equivalent
to saying that for any $ \phi\in L^1(E;m)$ with $\|\hat
G\phi\|_\infty<\infty$,
\begin{equation}
\label{eqi.144}
(u,\phi)=(f(\cdot,u),\hat G\phi)
+\int_E \hat G\phi \,d\mu+\int_E\hat G\phi\,d\nu
\end{equation}
(see \cite{KR:CM}). Note that (\ref{eqi.144}) is a generalization
of Stampacchia's definition by duality introduced in
\cite{Stampacchia2} for solutions of uniformly elliptic PDEs with
measure data. Another equivalent definition is given in
\cite{KR:NODEA}, where it is shown that (\ref{eq1.3}) is satisfied
if and only if $u$ is a renormalized solution to the first
equation of (\ref{eq1.2}), i.e. $u$ is quasi-continuous,
$f(\cdot,u)\in L^1(E;m)$, $T_k(u):=(u\wedge k)\vee (-k)$ belongs
to the extended Dirichlet space $D_e[\EE]$ (see Section \ref{sec2} for the
definition) and
\begin{equation}
\label{eq1.5} \EE(T_ku,v)=\int_Ef(\cdot,u)v\,dm+\int_E v\,d\mu
+\int_E v\,d\nu+\int_Ev\,d\nu_k
\end{equation}
for some sequence $\{\nu_k\}$ of bounded smooth measures on $E$
such that $\|\nu_k\|\rightarrow 0$ as $k\rightarrow\infty$,
where $\|\cdot\|$ stands for the total variation norm on the space
of signed Borel measures on $E$. The concept of renormalized
solutions to elliptic equations with measure data and local
operators of Leray-Lions type was introduced in \cite{DMOP}.

Our main result on existence and uniqueness of solutions of the
complementary system (\ref{eq1.2}) is first proved for one
reflecting barrier in Section \ref{sec3} and then for two barriers
in Section \ref{sec4}. It is worth mentioning that in both cases
we give necessary and sufficient conditions on barriers $h_{1},
h_{2}$ under which there exists a solution $u$ of (\ref{eq1.2})
with $f,\mu$ satisfying our assumptions. We also prove that $u$ is
an envelope  of supersolutions of some partial differential
equation related to (\ref{eq1.2}). More precisely, we show that
\begin{equation}
\label{eq1.4} u=\mbox{quasi-essinf$\{v\ge h_{1}$ a.e., $v$ is
a supersolution of PDE$(f+d\mu-d\nu^{-})$}\},
\end{equation}
where as before $\nu^-$ denotes the negative part of the reaction
measure.
A result similar to (\ref{eq1.4}) was proved in \cite{KR:JEE} for
evolution obstacle problem involving divergence form operator.

In case $\mu\in\MM_{0,b}$, $f(\cdot,0)\in L^1(E;m)$ and the
barriers satisfy some additional regularity condition we show that
$\nu\in\MM_{0,b}$. When combined with the regularity results
proved in \cite{KR:JFA,KR:CM} this  implies that for every $k\ge0$ the
truncation $T_k(u)$ of $u$ at the level $k$ belongs to the
extended Dirichlet space $D_e[\EE]$ and
\[
\EE(T_k(u),T_k(u))\le
2k(\|\mu\|+\|\nu\|+\|f(\cdot,0)\|_{L^1(E;m)}).
\]
Moreover, we show that if $u$ is a solution to (\ref{eq1.1})
and $\mu\in D'_e[\EE]$, where $D'_e[\EE]$ is the dual
of $D_e[\EE]$, and moreover, $f(\cdot,u)\in D'_e[\EE]$ and
there exists $v=R\lambda$  for some $\lambda\in D'_e[\EE]$ (in case of $h_2\equiv \infty$ it is enough to assume that $v\in D_e[\EE]$)  such that $h_1\le v\le
h_2$, then $u\in D_e[\EE]$, $\nu\in D'_e[\EE]$ and $(u,\nu)$ is the unique pair in  $D_e[\EE]\times D'_e[\EE]$ such that
\begin{equation}
\label{eq1.1234}
\EE(u,\eta)=\int_E f(\cdot,u)\eta\,dm
+\int_E\eta\,d\mu+\int_E\eta\,d\nu,\quad \eta\in D_e[\EE],
\end{equation}
\begin{equation}
\label{eq1.123}
 \int_E(u-h_1)\,d\nu^+= \int_E(h_2-u)\,d\nu^-=0,\quad h_1\le u\le h_2\quad\mbox{q.e.}
\end{equation}
This formulation of a solutions is equivalent to the variational
inequalities formulation i.e. finding $u\in D_e[\EE]$ such that $\psi_1\le u\le \psi_2,\, m$-a.e. and
\begin{equation}
\label{eq1.12345}
\EE(u,u-\eta)\le \int_E f(\cdot,u)(u-\eta)\,dm
+\int_E(u-\eta)\,d\mu,\quad \eta\in D_e[\EE],\, \psi_1\le\eta\le\psi_2.
\end{equation}
It is enough to put $(u-\eta)$ as test function in (\ref{eq1.1234}) and apply (\ref{eq1.123}).
Note here that in general it is not true that $L^2(E;m)$
is a subset of $D'_e[\EE]$.

In Section \ref{sec5}, we prove  a Lewy-Stampacchia type
inequality, which is known to be useful in the study of regularity
of solutions of (\ref{eq1.2}). If one of the barriers, say
$h_{1}$, is a difference of two natural potentials, then
\[
\nu^{+}\le\mathbf{1}_{\{u=h_{1}\}}(f(\cdot, h_{1})+\mu+Lh_{1})^{-}\cdot m.
\]
Note that even in the case of local operators there are only few
results of this type for two-sided obstacle problem (see
\cite{MM,MTV,RT}). We also prove some stability results which in particular implies
that probabilistic solutions to (\ref{eq1.1}) are pointwise limits of analytic solutions.

\section{Preliminaries}
\label{sec2}

For convenience of the reader and to fix notation, in this section
we provide some basic information on Dirichlet spaces and
associated Markov processes. For more details we refer the reader
to monographs \cite{CF,MR} (quasi-regular Dirichlet  forms) and \cite{Fukushima,O} (regular Dirichlet forms).

In the whole paper $E$ is a Lusin space and $m$ is a positive
$\sigma$-finite measure on the $\sigma$-field $\BB(E)$ of Borel subsets of $E$.

Let $D[\EE]$ be a dense linear subspace of $L^{2}(E,m)$ and let
$\EE:D[\EE]\times D[\EE]\rightarrow\BR$ be a bilinear form.

We say that $(\EE,D[\EE])$ is positive if $\EE(u,u)\ge 0$ for
$u\in D[\EE]$. A positive definite form $(\EE, D[\EE])$ is called
a coercive closed form if
\begin{enumerate}
\item [\rm{(a)}]$(\tilde{\EE},D[\EE])$ is a symmetric closed form
on $L^{2}(E;m)$, where $\tilde\EE$ denotes the symmetric part of
$\EE$, i.e. $\tilde{\EE}(u,v)=\frac12(\EE(u,v)+\EE(v,u))$, $u,v\in
D[\EE]$,
\item [\rm{(b)}] $(\EE,D[\EE])$ satisfies the weak sector
condition, i.e. there exists $K>0$ such that
\[
|\EE_1(u,v)|\le K\EE_1(u,u)^{1/2}\EE_1(v,v)^{1/2},\quad
u,v\in D[\EE].
\]
\end{enumerate}
Here and henceforth,
\[
\EE_{\alpha}(u,v)=\EE(u,v)+\alpha(u,v),\quad u,v\in D[\EE]
\]
for $\alpha>0$. A form  $(\EE,D[\EE])$ is said to satisfy the strong sector condition if there is $K>0$ such that
\[
|\EE(u,v)|\le K\EE(u,u)^{1/2}\EE(v,v)^{1/2},\quad
u,v\in D[\EE].
\]
Note that symmetric forms satisfy the strong sector condition  with $K=1$
by Schwarz's inequality.

We say that $(\EE,D[\EE])$ is a Dirichlet form if it is closed
coercive form and for all $u\in D[\EE],\, u^{+}\wedge 1\in D[\EE]$
and
\[
\EE(u+u^{+}\wedge 1,u-u^{+}\wedge1)\ge 0,
\quad \EE(u-u^{+}\wedge 1,u+u^{+}\wedge1)\ge 0.
\]

For a Dirichlet form $(\EE,D[\EE])$ there
exists a unique operator $(L,D(L))$ on $L^2(E;m)$ (sometimes
called Dirichlet operator) such that
\[
D(L)\subset D[\EE],\quad \EE(u,v)=(-Lu,v),\, u\in D(A),v\in
D[\EE].
\]
By $\{G_{\alpha}\}_{\alpha>0}$ (resp. $\{T_{t}\}_{t>0}$) we will
denote the strongly continuous contraction resolvent (resp.
semigroup) generated by $(L,D(L))$ (see \cite[Chapter I]{MR}).

Given $F\in\BB(E)$ we set $D[\EE]_{|F}=\{u\in D[\EE]: u=0\mbox{ on
}F^{c}\,\, m\mbox{-a.e.}\}$. An increasing sequence $\{F_{k}\}$ of
closed subsets of $E$ is called $\EE$-nest if $\bigcup_{k\ge 1}
D[\EE]_{|F_{k}}$ is dense in $D[\EE]$ with respect to the norm
$\tilde{\EE}_{1}^{1/2}$. A set $N$ is an $\EE$-exceptional set if
$N^{c}\subset\bigcap_{k\ge 1} F_{k}^{c}$ for some $\EE$-nest
$\{F_{k}\}$. We say that a property in $E$ holds q.e. if it holds
outside some exceptional set. By  \cite[Theorem
III.2.11]{MR} (see also \cite[Exercise III.2.3 ]{MR}), every Borel
$\EE$-exceptional set is of $m$ measure zero. Consequently, if
some property holds q.e., it holds $m$-a.e. For equivalent
definitions of $\EE$-nest and $\EE$-exceptional set, expressed in
terms of some capacity associated with $(\EE,D[\EE])$ we refer the
reader to \cite[Section III.2]{MR}.

For a given nest $\{F_{k}\}$ we set
\[
C(\{F_{k}\})=\{f:E\rightarrow\BR; f_{|F_{k}}\mbox{ is continuous
for every}\,k\ge 1\}.
\]
Similarly we define sets $L(\{F_{k}\}),\, U(\{F_{k}\})$ replacing in the above definition the
word "continuous" by lower semicontinuous (l.s.c. in abbreviation) and upper semicontinuous (u.s.c in abbreviation), respectively.
We say that a function $u$ on $E$ is $\EE$-quasi-continuous (resp. $\EE$-l.s.c., $\EE$-u.s.c.) if
there exists an $\EE$-nest $\{F_{k}\}$ such that $u\in
C(\{F_{k}\})$ (resp. $u\in L(\{F_{k}\})$, $u\in U(\{F_{k}\})$).

A Dirichlet form $(\EE,D[\EE])$ on $L^{2}(E;m)$ is called
quasi-regular if
\begin{enumerate}
\item [\rm{(a)}] there exists an $\EE$-nest $\{F_{k}\}$ consisting
of compact sets,
\item [\rm{(b)}] there exists an $\tilde{\EE}^{1/2}_{1}$-dense subset
of $D[\EE]$ whose elements have $\EE$-quasi-continuous
$m$-versions,
\item[\rm{(c)}]there exist a sequence $\{u_{n}\}\subset D[\EE]$
of $\EE$-quasi-continuous functions and an $\EE$-exceptional
set $N\subset E$ such that $\{u_{n}\}$ separates points of $E\setminus N$.
\end{enumerate}

Let $(\EE,D[\EE])$ be a quasi-regular Dirichlet form on
$L^{2}(E;m)$. Adjoin $\Delta$ as an extra point to $E$ and set $E_{\Delta}=E\cup\Delta$. It is known (see \cite[Chapter IV]{MR}) that there exists an $m$-tight special
standard process $\mathbf{M}=(\Omega,\FF,\{X_{t}\}_{t\ge 0}$,
$\{P_{x}\}_{x\in E_{\Delta}})$ with life time $\zeta$ properly associated
with the form $(\EE, D[\EE])$, i.e. for every $t>0$ and $f\in
\BB_{b}(E)\cap L^{2}(E;m)$,
\begin{equation}
\label{eq.sem}
T_{t}f(x)=E_xf(X_{t})
\end{equation}
for $m$-a.e. $x\in E$ and $x\mapsto E_xf(X_{t})$ is
$\EE$-quasi-continuous. Note that  $X_t=\Delta$, $t\ge\zeta$  and
that above and it what follows we admit the convention that each
function $f$ on $E$ is extended to $E_{\Delta}$ by putting
$f(\Delta)=0$. By $\mathcal T$ we denote the set of all stopping
times with respect to $\mathcal F$. In particular $\zeta\in
\mathcal T$.

We say that a positive measure $\mu$ on $\BB(E)$ is $\EE$-smooth
if $\mu(N)=0$ for every $\EE$-exceptional set $N\in\BB(E)$ and
there exists an $\EE$-nest $\{F_{k}\}$ of compact subsets of $E$
such that $\mu(F_{k})<\infty$ for $k\ge 1$. The set of all
$\EE$-smooth measures on $\BB(E)$ will be denoted by $S$. We
denote by $\MM_{0,b}$ the set of bounded Borel measures $\mu$ on
$E$ such that $|\mu|\in S$.

In the paper, we frequently use the notion of
additive functional (AF for short) of $\mathbf M$ (for the
definition see \cite[Section 5.1]{Fukushima}). We say that an AF
$A$ of $\mathbf M$ is positive (resp. continuous) if $A_t\ge0,\,
t\ge 0,\, P_x$-a.s. (resp. $t\rightarrow A_t$ is continuous on
$[0,\infty)$ $P_x$-a.s.) for q.e. $x\in E$. We say that a process
$A$ is a martingale AF of $\mathbf M$ if $A$ is an AF of $\mathbf
M$ and it is a martingale with respect to $\FF$ under the measure
$P_x$ for q.e. $x\in E$.

It is known (see \cite[Theorem VI.2.4]{MR}) that there is a
one-to-one correspondence between $\EE$-smooth measures and
positive continuous additive functionals (PCAFs) of $\mathbf{M}$.
This correspondence, called Revuz correspondence, can be expressed
as
\[
\lim_{t\searrow 0}E_{m}
\big(\frac{1}{t}\int_{0}^{t}f(X_{s})\,dA_{s}\big)
=\int_{E}f\,d\mu,\quad f\in\BB^{+}(E),
\]
where $E_m$ denotes the expectation with respect to the measure $P_m(\cdot)=\int_EP_x(\cdot)\,m(dx)$. For an $\EE$-smooth measure $\mu$ we denote by $A^{\mu}$ the  unique
PCAF of $\mathbf{M}$ associated with $\mu$. We also set for $\mu\in S$,
\[
R\mu(x)=E_x\int_0^\zeta \,dA^\mu_r,\quad x\in E.
\]
We say that a form $(\EE,D[\EE])$ is transient if the associated
semigroup $\{T_{t}\}_{t>0}$ is transient, i.e. $G\phi$  is finite
$m$-a.e. for every nonnegative $\phi\in L^{1}(E;m)$. Equivalently
(see \cite[Corollary 3.5.34]{J2}), the form is transient  if there
exists a strictly positive $g\in L^1(E;m)$ such that
(\ref{eq1.10}) is satisfied.

For a coercive closed form $(\EE,D[\EE])$ we define $D_e[\EE]$ as follows: $D_e[\EE]$ is
the family of all functions $u$ on $E$ for which there exists an
$\EE$-Cauchy sequence (i.e. Cauchy sequence with respect to the norm generated by the inner product $\tilde\EE$) $\{u_n\}\subset D[\EE]$ such that
$u_n\rightarrow u$ $m$-a.e. ($\{u_n\}$ is called the approximating
sequence for $u$). It is known that if $(\EE,D[\EE])$ is transient
then for each fixed $u\in D_e[\EE]$ the limit of $\{\EE(u_n,u_n)\}$ is
independent of the approximating sequence for $u$. We set
$\EE(u,u)=\lim_{n\rightarrow\infty}\EE(u_n,u_n)$. By \cite[Lemma 1.5.5]{Fukushima}, the pair $(\tilde\EE, D_e[\EE])$
is a Hilbert space. By \cite[Remark 2.2]{KR:arxiv},  each $u\in D_e[\EE]$ has an $m$-version which is quasi-continuous.
From now on for given $u\in D_e[\EE]$ we always consider its quasi-continuous $m$-version.

We denote by $\|\cdot\|_\EE$ the norm generated by $\tilde\EE$ and by $\|\cdot\|_{\EE'}$
the norm on its dual space. If $(\EE,D[\EE])$ is transient, then by \cite[Lemma 2.1]{KR:arxiv}, for every $\mu\in S$ there exists an $\EE$-nest $\{F_k\}$ such that
$\mathbf{1}_{F_k}\cdot\mu\in D_e'[\EE]$. If, in addition, $(\EE,D[\EE])$ satisfies the strong sector condition, then by
\cite[Lemma 2.4]{KR:CM}, if $\mu\in D'_e[\EE]$, then $u:=R\mu\in D_e[\EE]$ and
\begin{equation}
\label{eq2.eva}
\EE(u,\eta)=\int_E\eta\,d\mu,\quad \eta\in D_e[\EE].
\end{equation}

A nonnegative measurable function $u:E\rightarrow \BR$ is called
$\EE$-excessive if $T_{t}u\le u$ for  $t\ge 0$ $m$-a.e.
We say that $u$ is an $\EE$-natural potential if
there exists a positive  $\mu\in \mathbb M_0$ such that $u=R\mu$
q.e.  A function $f:E\rightarrow \BR$ is called
$\EE$-quasi-integrable ($f\in qL^{1}(E;m)$ in notation) if
$A^{|f|\cdot m}$ is a finite AF of $\mathbf{M}$. We say that
$f:E\rightarrow \BR$ is locally $\EE$-quasi-integrable if
$A^{|f|\cdot m}$ is an AF of $\mathbf{M}$.

In \cite{OP} the notion of quasi-integrability was considered in
the case of Laplace operator. Our notion of
quasi-integrability is more general (since it applies to wider
class of operators), but at the same time is stronger than the
notion introduced in \cite{OP} in the particular case of Laplace
operator. As a matter of fact, the quasi-integrability introduced
in \cite{OP} coincides with the local quasi-integrability
considered in the paper \cite{Kl:AMPA} devoted to elliptic systems
involving Laplace operator (see comments following \cite[Remark
2.3]{Kl:AMPA}). Note also that in the case of Laplace
operator the life-time $\zeta$ of the associated process is
predictable. Therefore the results of \cite{Kl:AMPA} suggest that
in the case of operators associated with a quasi-regular Dirichlet
form for which the life-time of the associated process is
predictable (e.g. regular Dirichlet form without killing part) the
main results of our paper hold true if in their assumptions we
replace quasi-integrability by local quasi-integrability.

In the paper, if there will be no ambiguity, we drop the
letter $\EE$ in the notation. For instance, instead of writing
$\EE$-quasi-continuous, $\EE$-smooth, etc. we simply write
quasi-continuous, smooth, etc. By $\rightarrow_P$ we denote the
convergence in probability $P$. $x^+=\max(x,0)$, $x^-=\max(-x,0)$.

\section{One-sided obstacle problem}
\label{sec3}

From now on,  $(\EE, D[\EE])$ is a
transient quasi-regular Dirichlet form satisfying the strong sector condition, $f:E\times\BR\rightarrow\BR$, $h,h_{1},h_{2}: E\rightarrow \BR$ are measurable functions and $\mu$ is a measure on $\BB(E)$ such that
$|\mu|\in S$.

Given $\mu\in S$ we define the  $0$-order potential operator by putting
\[
R\mu(x)=E_x\int^{\zeta}_0dA^{\mu}_t
\]
for q.e. $x\in E$.
In the important case where $\mu=f\cdot m$ for some $f\in L^1(E;m)$ the AF associated with $\mu$ has the form
$A^\mu_t=\int_0^tf(X_r)\,dr$, $t\ge0$ (see \cite[Theorem A.3.5]{CF} and remarks following it). Consequently, with our convention that $f(\Delta)=0$, in that case we have
\[
R\mu(x)=E_x\int_0^\infty f(X_t)\,dt
\]
for q.e $x\in E$.
From this and (\ref{eq.sem}) it follows that
\[
R\mu=Gf\quad m\mbox{-a.e.}
\]
The above relation may be easily extended to $f\in\BB^+(E)$ by approximation.

We will need the following hypotheses:
\begin{enumerate}
\item[(H1)] $y\mapsto f(x,y)$ is nonincreasing for every $x\in E$,
\item[(H2)]$y\mapsto f(x,y)$ is continuous for every $x\in E$,
\item[(H3)]$x\mapsto |f(x,y)|\in qL^{1}(E;m)$ for every $y\in\BR$,
\item[(H4)] $R|f(\cdot, 0)|+R|\mu|<\infty$ $m$-a.e.,
\item[(H5)] There exists $v:E\rightarrow\BR$ such that $v$ is a difference
of natural potentials and $m$-a.e.,
\[
v\ge h,\quad Rf^{-}(\cdot, v)<\infty,
\]
\item[(H6)] There exists $v:E\rightarrow\BR$ such that $v$ is a difference
of natural potentials and $m$-a.e.,
\[
h_{1}\le v\le h_{2},\quad R|f(\cdot, v)|<\infty.
\]
\end{enumerate}

\begin{remark}
 (i) Let $h\in \BB(E)$. If $C=\{u\in D(\EE): u\ge
h\}\neq\emptyset$, then there exists the smallest natural
potential $v\ge h$. This is a  consequence of the Lax-Milgram
theorem (see \cite[Proposition  III.1.5]{MR}). Therefore,  if
$C\neq\emptyset $ and $f^-(\cdot,v)\in L^1(E;m)$, then (H5) is
satisfied.
\\
 (ii) In practice, an effective criterion
ensuring (H6) is the following:
\begin{enumerate}
\item[(a)] $f^+(\cdot,h_1), f^-(\cdot,h_2)\in L^1(E;m)$,
\item[(b)] there exists $w\in D(L)$ and   $\varphi$ being a
difference of convex functions with $\varphi(0)=0$ such that
$h_1\le \varphi(w)\le
h_2$.
\end{enumerate}
By the Tanaka-Meyer formula  (see  \cite[Theorem IV.70]{Protter})
if (b) is satisfied, then  $\varphi(w)$ is a difference of natural
potentials.
\end{remark}

Let us define  the class $\mathbb M_0$ by (\ref{eq1.02}). In
\cite{KR:CM} it is shown that $\mathbb M_0$ can be equivalently
defined as
\begin{equation}
\label{eq3.01} \mathbb M_0=\{\mbox{$\mu:|\mu|\in S$,\,\,
$R|\mu|<\infty$ $m$-a.e.}\}.
\end{equation}
Note  also that from \cite[Corollary 1.3.6]{O} it follows
immediately that $\MM_{0,b}\subset\mathbb M_0$. So, we see that (H4)--(H6) are satisfied
in particular if $f(\cdot,0), f^-(\cdot,v)\in L^1(E;m), f(\cdot,v)\in L^1(E;m)$ and $\mu\in\MM_{0,b}$.
In general, the inclusion
is strict as the following examples show.

\begin{example}
Let $\alpha\in(0,2)$, $d\ge3$, and let $D\subset\BR^d$ be an open
bounded set with smooth boundary. Consider the form
$(\EE_D,D[\EE_D])$ associated with $\alpha$-Laplace operator
$\Delta_D^{\alpha/2}$ on $D$ with zero Dirichlet boundary
condition (see, e.g., \cite[Section 6.3]{KR:CM}). The form $\EE_D$ can be constructed as follows. We first consider the form $(\EE, D[\EE])$
associated with $\Delta^{\alpha/2}$ on $\BR^d$, i.e.
\[
\EE(u,v)=\int_{\BR^d}\hat u(x)\bar{\hat v}(x)\psi(x)\,dx,\, u,v\in D[\EE],
\]
where $\psi(x)=|x|^{\alpha/2}$ for $x\in\BR^d$ and   $\hat u,\hat v$  denote the  Fourier transforms of $u$ and $v$, and
\[
D[\EE]=\{w\in L^2(\BR^d):\int_{\BR^d}|\hat w(x)|^2\psi(x)\,dx<\infty\}
\]
(see \cite[Example 1.4.1]{Fukushima}). Next we set
\[
\EE_D(u,v)=\EE(u,v),\quad u,v\in D[\EE_D]:=\{w\in D[\EE]: w=0\mbox{ q.e. on }\, \BR^d\setminus D\},
\]
that is $(\EE_D,D[\EE_D])$ is the part of $(\EE,D[\EE])$ on $D$. By \cite[Theorems 4.4.3, 4.4.4]{Fukushima},   $(\EE_D,D[\EE_D])$ is again a regular symmetric transient Dirichlet form, so it generates
a Dirichlet operator which we denote by $\Delta_D^{\alpha/2}$. Note that from the definition of $D[\EE_D]$ it follows that in the case of  the nonlocal operator  $\Delta_D^{\alpha/2}$
zero boundary Dirichlet condition in fact means zero exterior condition.
By
\cite[Proposition 4.9]{Ku}, there exists constants $0<c_1<c_2$
depending only on $d,\alpha,D$ such that (\ref{eq1.11}) is satisfied
with $\delta(x)=\mbox{dist}(x,\partial D)$ and $G$ associated with $(\EE_D,D[\EE_D])$.
Therefore from (\ref{eq1.02}) immediately  follows that
$L^{1}(D;\delta^{\alpha/2}(x)\,dx) \subset\mathbb M_0$.
\end{example}

\begin{example}
\label{ex.ex} Let $(\EE,D[\EE])$ be a regular symmetric Dirichlet
form on $L^2(E;m)$ and let $\mu\in S$.  Consider the form
$(\EE^{\mu},D[\EE^{\mu}])$, the perturbation of $(\EE,D[\EE])$ by
$\mu$, which is defined by
\[
\EE^{\mu}(u,v)=\EE(u,v)+\int_Euv\,d\mu, \quad u,v\in
D[\EE^{\mu}]:=\{u\in D[\EE]: \int_E |u|^2\,d\mu<\infty\}.
\]
It is  known (see \cite[Section IV.4(c)]{MR} and \cite[Section 6.1]{Fukushima}) that
$(\EE^{\mu},D[\EE^{\mu}])$ is a quasi-regular Dirichlet form on
$L^2(E;m)$ and the $0$-order potential operator $R^{\mu}$
associated with $(\EE^{\mu},D[\EE^{\mu}])$ has the form
\[
R^{\mu}\nu(x)=E_{x}\int_{0}^{\infty} e^{-A^{\mu}_{t}}\,dA^{\nu}_t
\]
for $\nu\in S$ (here  $A^\mu, A^\nu$ are PCAFs of the process  $\mathbf{M}$   associated with ($\EE,D(\EE))$ in the Revuz correspondence with $\mu$ and $\nu$, respectively). In particular,
\[
R^{\mu}\mu(x)=E_{x}\int_{0}^{\infty} e^{-A^{\mu}_{t}}\,dA^{\mu}_{t}.
\]
The last integral is less than or equal to 1. Since by \cite[Lemma IV.4.5]{MR}  the measure $\mu$ is smooth with respect to the perturbed form $(\EE^{\mu},D[\EE^{\mu}])$, it follows from
(\ref{eq3.01}) that $\mu$ belongs to the class $\mathbb M_0(\EE^{\mu})$ defined for the form $(\EE^{\mu},D[\EE^{\mu}])$. This shows that even nowhere Radon measures may belong to the
class $\mathbb M_0$ (for construction of smooth  nowhere Radon measure see \cite[Section IV.4(c)]{MR}).
\end{example}

We denote by  $\mathfrak S^p_c$ the set of all quasi-continuous
functions on $E$ such that for q.e. $ x\in E$,
\[
E_{x}\sup_{t\ge0}\,|u(X_{t})|^{p}<\infty.
\]

\begin{definition}
We say that $u:E\rightarrow\BR$ is a solution of
PDE$(f+d\mu)$ if
\begin{enumerate}
\item[(a)] $u$ is quasi-continuous and
$f(\cdot,u)\cdot m\in\mathbb M_0$,
\item[(b)]for q.e. $x\in E$,
\[
u(x)=E_{x}\Big(\int_{0}^{\zeta}\,f(X_{t}, u(X_{t}))\,dt
+\int_{0}^{\zeta}\,dA_{t}^{\mu}\Big).
\]
\end{enumerate}
\end{definition}

\begin{definition}
We say that a pair $(u,\nu)$ is a solution of OP$(f+d\mu, h)$ if
\begin{enumerate}
\item[(a)] $u$ is quasi-continuous and
$\nu$, $f(\cdot,u)\cdot m\in\mathbb M_0$,
\item[(b)]for q.e. $x\in E$,
\begin{equation}
\label{eq3.08} u(x)=E_{x}\Big(\int_{0}^{\zeta}\,f(X_{t},
u(X_{t}))\,dt
+\int_{0}^{\zeta}\,dA_{t}^{\mu}+\int_{0}^{\zeta}\,dA_{t}^{\nu}\Big),
\end{equation}
\item[(c)] $u(x)\ge h(x)$ for $m$-a.e. $x\in E$,
\item[(d)] for q.e. $x\in E$,
\[
E_{x}\int_{0}^{\zeta}(u(X_{t})-h^{*}(X_{t}))\,dA_{t}^{\nu}=0
\]
for every quasi-continuous function $h^{*}$ on $E$ such that $h\le
h^{*}\le u$ $m$-a.e.
\end{enumerate}
\end{definition}

\begin{remark}
(i) By the Revuz duality,
condition (d) is equivalent to the following condition:
\[
\int_{E}(u-h^*)\,d\nu=0
\]
for every quasi-continuous function $h^*$ such that $h\le h^*\le
u$ $m$-a.e. Standard argument shows that in fact one can replace
$h^*$ by any quasi-u.s.c. $h^{**}$ such that $h\le h^{**}\le u$
$m$-a.e.
\smallskip\\
(ii) Let $\hat{h}$ be a quasi-u.s.c. regularization of $h$, i.e.
\[
\hat{h}=\mbox{quasi-essinf}\{\eta\ge h\,\,m\mbox{-a.e.}: \eta\mbox{ is
quasi-u.s.c.}\}.
\]
Then $(u,\nu)$ is a solution of OP$(f+d\mu,h)$ if and only if it is a
solution of OP$(f+d\mu,\hat{h})$. Indeed, if $(\hat{u},\hat{\nu})$
is a solution of OP$(f+d\mu,\hat{h})$ then of course $\hat{u}\ge
h$ $m$-a.e. Furthermore, for every quasi-u.s.c. $h^*$ such that $h\le
h^*\le\hat u$ $m$-a.e. we have
\[
\int_{E}(\hat{u}-h^*)\,d\hat{\nu}\le \int_{E}(\hat{u}-\hat{h})\,d\hat{\nu}=0
\]
since  $\hat{h}\le h^*$ q.e. Therefore $(\hat{u},\hat{\nu})$
is a solution of OP$(f+d\mu,h)$. Now assume that $(u,\nu)$ is a
solution of OP$(f+d\mu,h)$. Then $\hat{h}\le u$ q.e. since
$u$ is quasi-continuous, and
\[
\int_{E}(u-\hat{h})\,d\nu=0
\]
since $\hat{h}$ is quasi-u.s.c. and $h\le \hat{h}$ $m$-a.e. Thus
$(u,\nu)$ is a solution of OP$(f+d\mu,\hat{h})$. From the above it
follows that without loss of generality we can confine ourselves
to considering  quasi-u.s.c. barriers. Moreover, if $h$ is
quasi-u.s.c. then the minimality condition (d) reduces to
\[
\int_{E}(u-h)\,d\nu=0.
\]
\end{remark}

In the proof of Theorem \ref{th.m1} we will use the form
$(\EE^{\#},D[\EE^{\#}])$, which is described in detail in
\cite[Theorem VI.1.2]{MR}.  Here let us only mention that $E^{\#}$
is a local compactification of $E$ and  $(\EE^{\#},D[\EE^{\#}])$
is a regular Dirichlet form on $L^{2}(E^{\#};m^{\#})$, which is an
extension of the form $(\EE,D[\EE])$. Since $(\EE^{\#},
D[\EE^{\#}])$ is regular, one can associate with it  a Hunt
process
$\mathbf{M}^{\#}\equiv(\Omega^{\#},\FF^{\#},\{X^{\#}_{t}\}_{t\ge
0}, \{P^{\#}_{x}\}_{x\in E^{\#}_{\Delta}})$  with life time
$\zeta^{\#}$. The process $\mathbf{M}^{\#}$ being a Hunt process,
is a special standard process, and moreover, its trajectories have
left limits on $(0,+\infty)$. $\mathbf{M}^{\#}$ is a standard
extension of $\mathbf{M}$, i.e. $P_{x}=P^{\#}_{x}$,
$X_{t}=X^{\#}_{t},\, t\ge 0,\, P_{x}$-a.s. for every $x\in E$ and
$P^{\#}_{x}=\delta_{x}$, $X^{\#}_{t}=x$, $t\ge 0$, for every $x\in
E^{\#}\setminus E$. Given $u:E\rightarrow\BR$ we will denote by
$u^{\#}$ its extension to $E^{\#}$ defined as $u^{\#}(x)=u(x)$ for
$x\in E$ and $u^{\#}(x)=0$ for $x\in E^{\#}\setminus E$.

The above procedure  of regularization of
quasi-regular Dirichlet  form and associated Markov process is
called the transfer method in  \cite{MR}. In what follows, we use
this procedure without mentioning when we cite some results from
\cite{Fukushima} or other papers dealing with regular Dirichlet
forms (and not quasi-regular forms).

\begin{lemma}
\label{lm3.1} Suppose that $\mu\in\mathbb M_0$
and $u$ is a quasi-continuous function on $E$ such that
\begin{equation}
\label{eq3r.1} u(x)=E_x\int_0^\zeta\,dA^\mu_t
\end{equation}
for q.e. $x\in E$. Then there exists a martingale AF $M$ of
$\mathbf M$ such that for q.e. $x\in E$,
\begin{equation}
\label{eq3r.2}
u(X_t)=\int_t^\zeta\,dA^\mu_r-\int_t^\zeta\,dM_r,
\quad t\in [0,\zeta],\, P_x\mbox{a.s.}
\end{equation}

\end{lemma}
\begin{proof}
 By the transfer method, we may assume that
$\mathbf M$ is a  Hunt process. By \cite[Theorem
4.1.1]{Fukushima}, there exists a properly exceptional set
$N\subset E$ such that (\ref{eq3r.1}) holds for $x\in E\setminus
N$. Using  the Markov property and additivity of $A^\mu$ we
conclude  from (\ref{eq3r.1}) that
\begin{equation}
\label{eq3r.3} u(X_t)=E_{X_t}A^\mu_\zeta
=E_x(A^\mu_\zeta\circ\theta_t|\FF_t)
=E_x\Big(\int_0^\zeta\,dA^\mu_r\,|\,\FF_t\Big)-A^\mu_t,\quad t\ge
0,
\end{equation}
for every $x\in E\setminus N$. In the above equation, $\theta$ is
the  shift operator on $\Omega$, that is
$X_t(\theta_s\omega)=X_{s+t}(\omega),\, \omega\in\Omega,\, s,t\ge
0$. Set
\begin{equation}
\label{eq3r.4}
M_t=u(X_t)-u(X_0)+A^\mu_t,\quad t\ge 0.
\end{equation}
Clearly $M$ is an additive functional. By quasi-continuity of $u$
and \cite[Theorem 4.6.1]{Fukushima}, it is a c\`adl\`ag process.
By (\ref{eq3r.3}), $E_xM_t=0,\, t\ge 0,\, x\in E\setminus N$. Thus
$M$ is a martingale AF of $\mathbf M$. From (\ref{eq3r.4}) we get
(\ref{eq3r.2}).
\end{proof}

Now we will prove the main result of this
section. Besides the existence of a solution $(u,\nu)$ to
OP$(f+d\mu,h)$ we will show that $(u,\nu)$ can be approximated by
a solution $u_n$ to penalized  PDE (\ref{eq2.1}) with $\nu_n$
defined via $u_n$ and $h$. This approximation  is very important
in analysis of various properties of $u,\nu$ and in numerical
methods. We show the convergence of $u_n$ in the metric of the
space $\mathfrak S_c^q$, which implies the convergence of measures
$\nu_n$ to $\nu$ in the sense that  $A^{\nu_n}\rightarrow A^\nu$
in $\mathfrak S^q_c$ (clearly this convergence is stronger then
weak convergence since it preserves smoothness of measures). Note
here that in many applications the information about the measure
$\nu$ is crucial. As we have already mentioned in  Introduction,
$\nu$ can be interpreted as a local minimizer of the interacting
energy (\ref{eq1.00.2}). Moreover, in applications to mathematical
finance, the  AF $A^\nu$ generated by $\nu$  can be interpreted in
some models as the so-called early exercise premium (see
\cite{KR:MF}). As a by-product, we also get probabilistic
interpretation of solutions to OP$(f+d\mu,h)$. This result is a
basis for probabilistic numerical methods  (Monte Carlo methods)
and the optimal stopping theory, which  links  value functions of
type  (\ref{eq1.00.1}) with solutions  to OP$(f+d\mu,h)$.

Let $u$ be a real function on $E$. From now on,
\[
f_u(x):=f(x,u(x)),\quad x\in E.
\]

\begin{theorem}
\label{th.m1} Assume \mbox{\rm{(H1)--(H4)}}. Then there exists a
solution $(u,\nu)$ of \mbox{\rm{OP}}$(f+d\mu, h)$ if and only if
\mbox{\rm{(H5)}} is satisfied. Moreover, if \mbox{\rm(H5)} is
satisfied, then $u\in \mathfrak S^q_c$ for $q\in(0,1)$,
$u_{n}\rightarrow u$ in $\mathfrak S^q_c$ for $q\in(0,1)$ and
$u_{n}\nearrow u$ q.e., where $u_n$ is a unique solution of the
problem
\begin{equation}
\label{eq2.1}
-Lu_{n}=f(\cdot,u_{n})+\mu+\nu_{n}
\end{equation}
with $\nu_{n}=n(u_{n}-h)^{-}\cdot m$.
\end{theorem}
\begin{proof}
The necessity of (H5) follows from the fact that $u$ defined by
(\ref{eq3.08}) is a difference of natural potentials. To prove that (H5) is sufficient let us first note that
from \cite[Theorem 3.8]{KR:CM} (see also \cite[Theorem
4.7]{KR:JFA}) it follows that for each $n\in\BN$ there exists a
unique solution $u_{n}$ of (\ref{eq2.1}). Moreover, by
\cite[Proposition 4.9]{KR:JFA}, $u_{n}\le u_{n+1}$ q.e.
By (H5) there exists $\lambda\in\mathbb M_0$ such that
$-Lv=\lambda$ and $f^{-}(\cdot, v)\in\mathbb M_0$. Hence
\[
-Lv=\lambda+f_{v}+f_{v}^{-}-f_v^{+}.
\]
Let $\overline{v}$ be a solution of
\[
-L\overline{v}=\lambda^++f_{\overline{v}}+f^{-}_{v}+\mu^{+}.
\]
By \cite[Proposition 4.9]{KR:JFA}, $v\le\overline{v}$ q.e. Consequently,
$h\le\overline{v}$ q.e. From this we conclude that
\[
-L\overline{v}=\lambda^++f_{\overline{v}}+f^{-}_{v}
+\mu^{+}+n(\overline{v}-h)^{-}.
\]
By \cite[Proposition 4.9]{KR:JFA} again, for every $n\in\BN$,
\begin{equation}
\label{eq2.2} u_{n}\le\overline{v}\quad\mbox{q.e.}
\end{equation}
Set $u=\sup_{n\ge1}u_{n}$ and
\[
v_{n}(x)=-E_{x}\int_{0}^{\zeta}f(X_{t}, u_{n}(X_{t}))\,dt
-E_{x}\int_{0}^{\zeta}\,dA^{\mu}_t.
\]
Since $u_{n}\le u_{n+1}$ q.e., it follows from (H1) that $v_{n}\le
v_{n+1}$ q.e. For $n\in\BN$ set
\[
w_{n}(x)=u_{n}(x)+v_{n}(x).
\]
Then
\[
w_{n}(x)=E_{x}\,\int_{0}^{\zeta}\,dA^{\nu_{n}}_{t}.
\]
From this we see that $w_{n}$ is a
natural potential. In particular, $w_{n}$ is an excessive function.
Therefore  $w$ defined as
\[
w(x)=\sup_{n\ge 1}w_{n}(x) \quad \mbox{for q.e. } x\in E
\]
is excessive too (see \cite[Proposition 1.2.1]{BB}), and hence quasi-continuous (see  \cite[Theorem A.2.7]{Fukushima} and  \cite[Theorem 4.6.1]{Fukushima}). By (\ref{eq2.2}), (H1), (H2) and the Lebesgue dominated convergence theorem, we have
\begin{equation}
\label{eq2.3} v_{n}(x)\rightarrow -E_{x}\int_{0}^{\zeta}f(X_{t},
u(X_{t}))\,dt -E_{x}\int_{0}^{\zeta}\,dA_{t}^{\mu}.
\end{equation}
Hence
\[
w(x)=u(x)-E_{x}\int_{0}^{\zeta}f(X_{t}, u(X_{t}))\,dt
-E_{x}\int_{0}^{\zeta}\,dA_{t}^{\mu}
\]
for q.e. $x\in E$. From the above equation, (\ref{eq2.2}), quasi-continuity of $w$ and \cite[Theorem VI.4.22]{BG}
we conclude that $w$ is a natural potential.  Therefore
there exists a smooth measure $\nu$ such that for q.e. $x\in E$,
\[
w(x)=E_{x} \int_{0}^{\zeta} dA_{t}^{\nu}.
\]
Hence
\[
u(x)=E_{x}\int_{0}^{\zeta}f(X_{t}, u(X_{t}))\,dt
+E_{x}\int_{0}^{\zeta}\,dA_{t}^{\mu}
+E_{x}\int_{0}^{\zeta}\,dA_{t}^{\nu}
\]
for q.e. $x\in E$. By Lemma \ref{lm3.1} there exists a martingale AF $M$ of $\mathbf{M}$
such that
\[
u(X_{t})=\int_{t}^{\zeta}f_{u}(X_{r})\,dr
+\int_{t}^{\zeta}dA_{r}^{\mu}+\int_{t}^{\zeta}dA_{r}^{\nu} +\int_{t}^{\zeta}dM_{r},
\quad 0\le t\le\zeta,\quad P_x\mbox{-a.s.}
\]
for q.e. $x\in E$. Since $u_{n}, u$ are
quasi-continuous and we know that $u_{n}\rightarrow u$ and
$u_{n}\le u_{n+1}$ q.e., we see that $u^{\#}, u_{n}^{\#}$ are
$\EE^{\#}$-quasi-continuous, $u^{\#}_{n}\rightarrow u^{\#}$ and
$u^{\#}_{n}\le u_{n+1}^{\#}$, $\EE^{\#}$-q.e.  Therefore by
\cite[Theorem IV.5.29]{MR}, $u^{\#}_{n}(X^{\#}_{t})\rightarrow
u^{\#}(X^{\#}_{t})$, $t\ge 0$, and
$u^{\#}_{n}(X^{\#}_{t-})\rightarrow u^{\#}(X^{\#}_{t-})$, $t\ge
0$, $P^{\#}_{x}$-a.s. for $\EE^{\#}$-q.e. $x\in E^{\#}$. By
\cite[Proposition V.2.28]{MR} (see also \cite[Proposition
V.2.12]{MR}),
$u^{\#}_{n}(X^{\#}_{t-})=(u^{\#}_{n}(X^{\#}_{t}))_{-}$ and
$u^{\#}(X^{\#}_{t-})=(u^{\#}(X^{\#}_{t}))_{-}$ for $t\ge 0$.
Therefore by Dini's theorem, for every $T>0$,
\[
\sup_{t\le T}|u^{\#}_{n}(X^{\#}_{t})-u^{\#}(X^{\#}_{t})|
\rightarrow_{P^{\#}_{x}} 0
\]
for $\EE^{\#}$-q.e. $x\in E$, which implies that
\begin{equation}
\label{eq2.4} \sup_{t\le T}
|u_{n}(X_{t})-u(X_{t})|\rightarrow_{P_{x}} 0
\end{equation}
for $\EE$-q.e. $x\in E$. Since the finite variation parts of the
semimatringales $u_{0}(X)$ and $u(X)$ are continuous,  $u_{0}(X), u(X)$
are special semimartingales (see \cite[Theorem III.34]{Protter}). Therefore there
exists an increasing  sequence $\{\tau_{k}\}\subset \mathcal{T}$
such that $\tau_k\nearrow \infty$, and
\[
E_{x}\sup_{t\le \tau_{k}}|u(X_{t})|
+E_{x}\sup_{t\le \tau_{k}}|u_{0}(X_{t})|<\infty,\quad k\ge 1.
\]
Since $u_{0}\le u_{n}\le u$ for $n\ge 1$, (H1) implies that for
q.e. $x\in E$,
\begin{align}
\label{eq2.41}
\nonumber E_{x}\int_{0}^{\zeta}dA_{t}^{\nu_{n}}
& \le E_{x}\sup_{t\le \tau_{k}}|u(X_{t})| +
E_{x}\sup_{t\le \tau_{k}}|u_{0}(X_{t})|+
E_{x}\int_{0}^{\tau_{k}} |f_{u}(X_{t})|\,dt \\&
\quad+E_{x}\int_{0}^{\tau_{k}}|f_{u_{0}}(X_{t})|\,dt
+E_{x}\int_{0}^{\tau_{k}}dA_{t}^{|\mu|}.
\end{align}
This when combined with (\ref{eq2.4}) implies that for every
$T>0$,
\[
\left[u_{n}(X)-u(X)\right]_{T}
=\left[M^{n}-M\right]_{T}\rightarrow_{P_{x}}0
\]
(see \cite[Theorem 1.8]{Jacod}), which is equivalent (since $\sup_{t\le\tau_k}|\Delta M_t^n|$ is uniformly integrable with respect to $n$) to
\begin{equation}
\label{eq2.5}
\sup_{t\le T}|M_{t}^{n}-M_{t}|\rightarrow_{P_{x}} 0.
\end{equation}
By (\ref{eq2.2}), (H1), (H2) and the Lebesgue dominated convergence theorem, we get
\begin{equation}
\label{eq3.07}
E_{x}\int_{0}^{\zeta}|f_{u_{n}}(X_{t})-f_{u}(X_{t})|\,dt
\rightarrow0.
\end{equation}
From (\ref{eq2.4}), (\ref{eq2.5}) and (\ref{eq3.07}) it follows
that for every $T>0$,
\begin{equation}
\label{eq2.6}
\sup_{t\le T}|\int_{0}^{t}dA_{t}^{\nu_{n}}-\int_{0}^{t}dA_{t}^{\nu}|
\rightarrow_{P_{x}} 0
\end{equation}
for q.e. $x\in E$. Observe that by (\ref{eq2.41}),
\[
E_{x}\int_{0}^{\tau_{k}}(u_{n}(X_{t})-h(X_{t}))^{-}\,dt \rightarrow0
\]
for q.e. $x\in E$, which when combined with (\ref{eq2.4}) implies
that $u\ge h$ $m$-a.e. Finally, let $h^{*}$ be a quasi-continuous
function such that $h\le h^{*}\le u$ $m$-a.e. Then by
(\ref{eq2.4}) and (\ref{eq2.6}), for every $T>0$ we have
\[
\int_{0}^{T}(u_{n}(X_{t})-h^{*}(X_{t}))^{+}\,dA_{t}^{\nu_{n}}
\rightarrow_{P_{x}} \int_{0}^{T}(u(X_{t})-h^{*}(X_{t}))^{+}\,dA_{t}^{\nu}.
\]
On the other hand,
\begin{align*}
\int_{0}^{T}(u_{n}(X_{t})&-h^{*}(X_{t}))^{+}\,dA_{t}^{\nu_{n}}
\\&= n \int_{0}^{T}(u_{n}(X_{t})-h^{*}(X_{t}))^{+}
(u_{n}(X_{t})-h(X_{t}))^{-}\,dt \le 0,
\end{align*}
which implies that
\[
\int_{0}^{T}(u(X_{t})-h^{*}(X_{t}))\,dA^{\nu}_t=0 \quad
P_{x}\mbox{-a.s.}
\]
since $h^{*}\le u$ q.e. Therefore, $(u,\nu)$ is a solution to OP$(f+d\mu,h)$.
By \cite[Theorem 3.8]{KR:CM}, $u_n, \bar v\in\mathfrak S^q_c$ for every $q\in (0,1)$.
From this, (\ref{eq2.2}) and (\ref{eq2.4}),  we conclude that $u\in \mathfrak S^q_c,\, q\in (0,1)$, and   $u_{n}\rightarrow u$ in
$\mathfrak S^q_c$ for $q\in (0,1)$.
This completes the proof.
\end{proof}
\medskip

\begin{corollary}
\label{wn.wn.wn.wn} Assume \mbox{\rm(H1)--(H5)}
and retain the notation from Theorem \ref{th.m1}. Then for every
$q\in(0,1)$,
\[
\lim_{n\rightarrow\infty}E_{x}\sup_{t\ge
0}|A_{t}^{\nu_{n}}-A_{t}^{\nu}|^{q} \rightarrow 0
\]
for q.e. $x\in E$.
\end{corollary}
\begin{proof}
 Follows from (\ref{eq2.6}) and (\ref{eq2.41}).
\end{proof}

In what follows we denote by $\|\cdot\|$ the total variation norm
on the space of signed Borel measures on $E$.

\begin{proposition}
\label{stw.r1} Assume \mbox{\rm{(H1)--(H5)}}. Let $(u,\nu)$ be a
solution of \mbox{\rm{OP}}$(f+d\mu,h)$. Then
\[
\|\nu\|\le 2(\|\mu\|+\|f_0\|+\|\lambda^+\|+\|f_v^-\|)
\]
with $\lambda=-Lv$, where $v$ is the function from condition
\mbox{\rm(H5)}.
\end{proposition}
\begin{proof}
Let $\bar v$ be as in the proof of Theorem \ref{th.m1}.
By (\ref{eq2.2}),
\begin{align*}
E_x\int_0^{\zeta}dA^{\nu}_t&\le
E_x\int_0^{\zeta}dA^{\mu^-}_t+E_x\int_0^{\zeta}f^-_u(X_t)\,dt
+E_x\int_0^{\zeta}dA^{\lambda^+}_t\\
&\quad+E_x\int_0^{\zeta}f^+_{\bar{v}}(X_t)\,dt
+E_x\int_0^{\zeta}f^-_v(X_t)\,dt.
\end{align*}
By \cite[Lemma 2.6]{KR:CM} (see also \cite[Lemma 5.4]{KR:JFA}),
\[
\|\nu\|\le \|\mu^-\|+\|f_u^-\|+\|\lambda^+\|
+\|f^+_{\bar{v}}\|+\|f^-_v\|.
\]
By (H1) and (\ref{eq2.2}), $f^-_{u}\le f^-_{\bar{v}}$. Therefore
\[
\|\nu\|\le \|\mu^-\|+\|\lambda^+\|+\|f_{\bar{v}}\|+\|f^-_v\|.
\]
Since by \cite[Proposition 3.10]{KR:CM}, $\|f_{\bar{v}}\|\le
\|\lambda^+\|+\|f_v^-\|+\|\mu^+\|+\|f_0\|$, the desired estimate follows.
\end{proof}
\medskip

For  $k\ge 0$ we define the truncation  operator
$T_k:\mathbb{R}\rightarrow\mathbb{R}$ as
\[
T_k(y)=\min\{\max\{-k,y\},k\},\quad y\in\mathbb{R}.
\]

\begin{proposition}
Assume \mbox{\rm{(H1)--(H5)}}. Let $(u,\nu)$ be a solution  of
\mbox{\rm{OP}}$(f+d\mu,h)$. If $f^-_{v},\mu,\lambda^+,f_0\in
\mathcal{M}_{0,b}$ then $\nu\in\mathcal{M}_{0,b}$,
$T_{k}(u)\in D_e[\EE]$ for every $k\ge0$, and
\begin{equation}
\label{eq3.13}
\EE(T_k(u),T_k(u))\le2k(\|\mu\|+\|\nu\|+\|f_0\|),\quad k\ge 0.
\end{equation}
\end{proposition}
\begin{proof}
Follows from Proposition \ref{stw.r1} and \cite[Proposition 3.10,
Theorem 4.2]{KR:CM}.
\end{proof}
\medskip

The uniqueness of solutions of the obstacle problem follows from
the following comparison result, in which we assume that $f_{1},
f_{2}: E\times \mathbb{R}\rightarrow \mathbb{R}$, $h_{1},
h_{2}:E\rightarrow \BR$ are measurable and
$\mu_{1},\mu_{2}\in\mathbb M_0$.

\begin{proposition}
\label{stw3.33} Assume that $(u_{i},\nu_{i})$, $i=1,2$, is a
solution of \mbox{\rm{OP}}$(f_{i}+d\mu_{i}, h_{i})$. If
\[
d\mu_{1}\le d\mu_{2},\quad h_{1}\le h_{2}\,\,m\mbox{\rm-a.e.}
\]
and either
\[
\mbox{$f_{1}$\ satisfies {\rm(H1)} and  $f_{1}(\cdot, u_{2})\le
f_{2}(\cdot,u_{2})$ $m$\mbox{\rm-a.e.}}
\]
or
\[
\mbox{$f_{2}$ satisfies {\rm(H1)} and $f_{1}(\cdot,u_{1})\le
f_{2}(\cdot,u_{1})$ $m$\mbox{\rm-a.e.}},
\]
then $ u_{1}\le u_{2}$ q.e. Moreover, if $h_{1}=h_{2}$ and $f_1, f_2$ satisfy \mbox{\rm(H1)}, then
$d\nu_{1}\ge d\nu_{2}$.
\end{proposition}
\begin{proof}
Suppose that $f_{1}$ satisfies (H1) and $f_{1}(\cdot, u_{2})\le
f_{2}(\cdot, u_{2})$ $m$-a.e. Since the  Revuz correspondence is
one-to-one, we have
\begin{equation*}
\int_{0}^{t}f_{1}(X_{r}, u_{2}(X_{r}))\,dr \le
\int_{0}^{t}f_{2}(X_{r}, u_{2}(X_{r}))\,dr, \quad
\int_{0}^{t}dA_{r}^{\mu_{1}}\le \int_{0}^{t}dA_{r}^{\mu_{2}},
\,\, t\ge0.
\end{equation*}
By the definition of a solution to the obstacle problem and  Lemma \ref{lm3.1}
there exist  martingale AFs $M^1,M^2$ of $\mathbf{M}$
such that
\[
u_i(X_{t})=\int_{t}^{\zeta}f_{i}(X_r,u_i(X_r))\,dr
+\int_{t}^{\zeta}dA_{r}^{\mu_i}+\int_{t}^{\zeta}dA_{r}^{\nu_i} +\int_{t}^{\zeta}dM^i_{r},
\quad 0\le t\le\zeta,
\]
$P_x\mbox{-a.s.},\, i=1,2$ for q.e. $x\in E$.
By the Tanaka-Meyer formula (see, e.g., \cite[Theorem IV.70]{Protter}),
for every $\tau\in\mathcal{T}$ we
have
\begin{align*}
(u_{1}-u_{2})^{+}(X_{t}) &\le (u_{1}-u_{2})^{+}(X_{\tau}) \\
&\quad+ \int_{t}^{\tau}(f_{1}(X_{r}, u_{1}(X_{r}))- f_{2}(X_{r},
u_{2}(X_{r})))\mathbf{1}_{\{u_{1}>u_{2}\}}(X_{r})\,dr \\& \quad +
\int_{t}^{\tau}\mathbf{1}_{\{u_{1}>u_{2}\}}(X_{r})\,
d(A_{r}^{\mu_{1}}-A_{r}^{\mu_{2}}) +
\int_{t}^{\tau}\mathbf{1}_{\{u_{1}>u_{2}\}}(X_{r})\,dA_{r}^{\nu_{1}}\\
&\quad-
\int_{t}^{\tau}\mathbf{1}_{\{u_{1}>u_{2}\}}(X_{r})\,dA_{r}^{\nu_{2}}
-\int_{t}^{\tau}\mathbf{1}_{\{u_{1}>u_{2}\}}(X_{r-})\,
d(M_{r}^{1}-M_{r}^{2})\\
&=:\sum_{i=1}^{6} I_{i}(t, \tau).
\end{align*}
Observe that $I_{2}(t, \tau)\le 0$ by the assumptions on $f_1,f_2$. Since
$h_{1}\le u_{1} \wedge u_{2}\le u_{1}$,
\[
I_{4}(t, \tau)= \int_{t}^{\tau}(u_{1}-u_{2})^{-1}\cdot
(u_{1}-u_{1}\wedge u_{2})(X_{r})
\mathbf{1}_{\{u_{1}>u_{2}\}}(X_{r})\,dA_{r}^{\nu_{1}}=0.
\]
It is also clear that $I_{3}(t,\tau)\le 0$ and
$I_{5}(t,\tau)\le0$. Let $\{\tau_{k}\}\subset \mathcal{T}$ be a
fundamental sequence for the martingale $M^{1}-M^{2}$. Then by the
above estimates,
\[
E_{x}(u_{1}-u_{2})^{+}(X_{t\wedge\tau})\le
E_{x}(u_{1}-u_{2})^{+}(X_{\tau})
\]
for q.e. $x\in E$. From this and the fact that $u_{1}, u_{2}$ are
differences of natural potentials we conclude that $u_{1}\le
u_{2}$ q.e. Now assume that $h_{1}=h_{2}$. By Corollary \ref{wn.wn.wn.wn},
for every $T>0$,
\[
\sup_{t\le T}|A_{t}^{\nu_{n}^{1}}-A_{t}^{\nu_{1}}| + \sup_{t\le
T}|A_{t}^{\nu_{n}^{2}}-A_{t}^{\nu_{2}}| \rightarrow_{P_{x}} 0
\]
for q.e. $x\in E$, where $u_{n}^{i}$ is a solution of
\[
-Lu_{n}^{i}=f_{i}(x, u_{n}^{i}) +\mu_{i} + n(u_{n}^{i}-h_{1})^{-}
\]
and $\nu_{n}^{i}=n(u_{n}^{i}-h_{1})^{-}\cdot m$. By
\cite[Proposition 4.9]{KR:JFA}, $u_{n}^{1}\le u_{n}^{2}$ q.e.,
which implies the second assertion of the proposition.
\end{proof}

\begin{corollary}
Under \mbox{\rm{(H1)}} there exists at most one solution of
\mbox{\rm{OP}}$(f+d\mu, h)$.
\end{corollary}

In the case where $L$ is a uniformly elliptic divergence form
operator with zero Dirichlet boundary conditions the existence and
uniqueness of a solution $(u,\nu)$ to the problem (\ref{eq1.1})
(in the sense of the definition of the present paper) was proved
in \cite{RS}. In \cite{RS} it is assumed that $h$ is
quasi-continuous, $\mu\in\MM_{0,b}$ and $f$ satisfies (H1), (H2)
and  slightly stronger than (H3)--(H5) integrability conditions.
Note also that in the special case considered in \cite{RS}, $u$ is
an entropy solution of (\ref{eq1.1}).

\begin{definition}
We say that $v$ is a supersolution of PDE$(f+d\mu)$ if there
exists a positive $\lambda\in\mathbb M_0$ such that $v$ is a
solution of PDE$(f+d\mu+d\lambda)$.
\end{definition}

\begin{proposition}
Assume \mbox{\rm{(H1)--(H4)}}. Let $u$ be a solution of
\mbox{\rm{OP}}$(f+d\mu, h)$. Then
\[
u= \mbox{\rm quasi-essinf}\{v\ge h\,\, m\mbox{-a.e.}: v\mbox{ is a supersolution of \rm PDE}(f+d\mu)\}.
\]
\end{proposition}
\begin{proof}
Let $v$ be a supersolution of PDE$(f+d\mu)$ and $v\ge h$ $m$-a.e.
Then
\[
-Lv=f(\cdot, v)+\mu+\lambda+n(v-h)^{-}.
\]
By  $u_n$ denote the  solution of
\[
-Lu_{n}=f(\cdot, u_{n})+\mu+n(u_{n}-h)^{-}.
\]
By \cite[Proposition 4.9]{KR:JFA}, $u_{n}\le v$. Since we know that  $u_{n}\nearrow u$ q.e., the desired assertion follows.
\end{proof}

\begin{proposition}
\label{stw3.var}
Let $(u,\nu)$ be a solution to \mbox{\rm OP$(d\mu,h)$}. Assume that $\mu\in\ D_e'[\EE]$, and there exists $v\in D_e[\EE]$ such that $v\ge h$. Then $u\in D_e[\EE]$, $\nu\in D'_e[\EE]$ and
$(u,\nu)$ is the unique pair in $D_e[\EE]\times D_e'[\EE]$ such that
\begin{equation}
\label{eq3.var}
\EE(u,\eta)=\int_E\eta\,d\mu+\int_E\eta\,d\nu,\quad \eta\in D_e[\EE],\qquad u\ge h\,\,\mbox{a.e.}
\end{equation}
and
\[
\int_E(u-\eta)\,d\nu\le 0,\quad\eta\in D_e[\EE],\, \eta\ge h ,\,\mbox{a.e.}
\]
Moreover,
\begin{equation}
\label{eq3.ine}
\|\nu\|_{\EE'}\le 3(\|v\|_\EE+\|\mu\|_{\EE'}).
\end{equation}
\end{proposition}
\begin{proof}
By Theorem \ref{th.m1}, $u_n\nearrow u$, where
\begin{equation}
\label{eq3.ppp}
-Lu_n=\mu+\nu_n,\qquad \nu_n=n(u_n-h)^-\cdot m.
\end{equation}
By the definition of a solution to (\ref{eq3.ppp}),
\[
u_n=R\mu+R\nu_n.
\]
Let $\{F_k\}$ be an $\EE$-nest such that $\nu^k_n=\mathbf{1}_{F_k}\cdot\nu_n\in D'_e[\EE]$, and let
\begin{equation}
\label{eq3.mid}
u_n^k=R\mu+R\nu_n^k.
\end{equation}
By (\ref{eq2.eva}), $u_n^k\in D_e[\EE]$ and
\begin{equation}
\label{eq3.varn}
\EE(u_n^k,\eta)=\int_E\eta\,d\mu+\int_E\eta\,d\nu^k_n,\quad \eta\in D_e[\EE].
\end{equation}
Setting $\eta=u_n^k-v$ and using the fact that $\int_E(u_n^k-v)\,d\nu_n^k\le 0$ we easily get
\begin{equation}
\label{eq3.17}
\|u_n^k\|_\EE\le 2(\|v\|_\EE+\|\mu\|_{\EE'}).
\end{equation}
Let $\eta\in D_e[\EE]$ be a positive function. Then
\begin{equation}
\label{eq3.18}
\int_E\eta\,d\nu_n^k=\EE(u_n^k,\eta)-\int_E\eta\,d\mu\le \|u_n^k\|_\EE\|\eta\|_\EE+\|\eta\|_\EE\|\mu\|_{\EE'}.
\end{equation}
From (\ref{eq3.mid}) it is clear that $u^k_n\rightarrow u_n$  q.e.  as $k\rightarrow\infty$. Since  $(\EE, D_e[\EE])$ is a Hilbert space, it follows from this and (\ref{eq3.17}) that $u^k_n\rightarrow u_n$ weakly in $(\EE, D_e[\EE])$ as $k\rightarrow\infty$. On the other hand,
$\int_E\eta\,d\nu^k_n\rightarrow\int_{\bigcup^{\infty}_{k=1}F_k} \eta\,d\nu_n=\int_E\eta\,d\nu_n$, the equality being a consequence of the fact that $E\setminus\bigcup^{\infty}_{k=1}F_k$ is  $\EE$-exceptional.
 Therefore letting $k\rightarrow\infty$ in (\ref{eq3.varn}) shows that
\begin{equation}
\label{eq3.19}
\EE(u_n,\eta)=\int_E\eta\,d\mu+\int_E\eta\,d\nu_n,\quad \eta\in D_e[\EE].
\end{equation}
Furthermore, by (\ref{eq3.17}) and (\ref{eq3.18}),
\begin{equation}
\label{eq3.20}
\|u_n\|_\EE\le 2(\|v\|_\EE+\|\mu\|_{\EE'}),\qquad  \int_E\eta\,d\nu_n\le \|u_n\|_\EE\|\eta\|_\EE+\|\eta\|_\EE\|\mu\|_{\EE'}.
\end{equation}
Similarly, since $u_n\nearrow u$, it follows from the first inequality in (\ref{eq3.20})  that  $\{u_n\}$ is weakly convergent in
$(\EE, D_e[\EE])$ to $u\in D_e[\EE]$. From (\ref{eq3.20}) it also follows that, up to a subsequence,  $\{\nu_n\}$ is weakly convergent in $(\EE, D'_e[\EE])$ to $\tilde\nu\in D'_e[\EE]$. Letting $n\rightarrow\infty$ in (\ref{eq3.19}) we obtain  the variational equality in (\ref{eq3.var}) with $\nu$ replaced by $\tilde\nu$. By virtue of  (\ref{eq2.eva}) this implies that
\[
u=R\mu+R\tilde\nu\quad\mbox{q.e.},
\]
so $R\nu=R\tilde\nu$, q.e.,  which forces $\tilde\nu=\nu$.
 By this and (\ref{eq3.20}),  $\nu$ satisfies (\ref{eq3.ine}). The
other properties of $(u,\nu)$ formulated in (\ref{eq3.var}) follow from the definition of a solution of OP$(d\mu,h)$.
\end{proof}

\section{Two-sided obstacle problem}
\label{sec4}

\begin{definition}
We say that a pair $(u,\nu)$ is a solution of
\mbox{\rm{OP}}$(f+d\mu, h_{1}, h_{2})$ if
\begin{enumerate}
\item[(a)] $u$ is quasi-continuous and
$\nu \in \mathbb M_0$, $f(\cdot,u)\cdot m\in\mathbb M_0$,
\item[(b)]for q.e. $x\in E$,
\[
u(x)=E_{x}\Big(\int_{0}^{\zeta}f(X_{t}, u(X_{t}))\,dt
+\int_{0}^{\zeta}dA_{t}^{\mu}+\int_{0}^{\zeta}dA_{t}^{\nu}\Big).
\]
\item[(c)]$h_{1}(x)\le u(x)\le h_{2}(x)$ for  $m$-a.e. $x\in E$,
\item[(d)]for q.e. $x\in E$,
\[
E_{x}\int_{0}^{\zeta}(u(X_{t})-h_{1}^{*}(X_{t}))\,dA_{t}^{\nu^{+}}=
E_{x}\int_{0}^{\zeta}(h_{2}^{*}(X_{t})-u(X_{t}))\,dA_{t}^{\nu^{-}}=0
\]
for any quasi-continuous functions $h_{1}^{*}, h_{2}^{*}$ on $E$
such that $h_{1}\le h_{1}^{*}\le u\le h_{2}^{*}\le h_{2}$ $m$-a.e.
\end{enumerate}
\end{definition}

\begin{proposition}
\label{stw4.44} Let $(u_{i},\nu_{i})$, $i=1,2$, be a solution of
\mbox{\rm{OP}}$(f_{i}+d\mu_{i}, h_{1}^{i}, h_{2}^{i})$. Assume
that
\[
d\mu_{1}\le d\mu_{2},\quad h_{1}^{1}\le h_{1}^{2},\quad
h_{2}^{1}\le h_{2}^{2}\,\, m\mbox{\rm-a.e.}
\]
and either
\[
f_{1}\mbox{ satisfies  \mbox{\rm(H1)} and $f_{1}(\cdot,u_{2})\le
f_{2}(\cdot,u_{2})$ $m$\mbox{\rm-a.e.}}
\]
or
\[
\mbox{$f_{2}$ satisfies \mbox{\rm(H1)} and $f_{1}(\cdot,u_{1})\le
f_{2}(\cdot,u_{1})$ $m$\mbox{\rm-a.e.}}
\]
Then $u_{1}(x)\le u_2(x)$ for q.e. $x\in E$.
\end{proposition}
\begin{proof}
Since the Revuz correspondence is one-to-one,
\[
\int_{0}^{t}f_{1}(X_{r}, u_{2}(X_{r}))\,dr \le
\int_{0}^{t}f_{2}(X_{r}, u_{2}(X_{r}))\,dr, \quad
\int_{0}^{t}dA_{r}^{\mu_{1}} \le
\int_{0}^{t}dA_{r}^{\mu_{2}}, \,\, t\ge 0.
\]
By the definition of solution to the obstacle problem and  Lemma \ref{lm3.1},
there exist  martingale AFs $M^1,M^2$ of $\mathbf{M}$
such that
\[
u_i(X_{t})=\int_{t}^{\zeta}f_{i}(X_r,u_i(X_r))\,dr
+\int_{t}^{\zeta}dA_{r}^{\mu_i}+\int_{t}^{\zeta}dA_{r}^{\nu_i} +\int_{t}^{\zeta}dM^i_{r},
\quad 0\le t\le\zeta,
\]
$ P_x\mbox{-a.s.},\,i=1,2$ for q.e. $x\in E$.
By the Tanaka-Meyer formula (see \cite[Theorem IV.70]{Protter}), for
every $\tau\in\mathcal{T}$,
\begin{align*}
(u_{1}-u_{2})^{+}(X_{r})&\le (u_{1}-u_{2})^{+}(X_{\tau})\\
&\quad+\int_{t}^{\tau}\mathbf{1}_{\{u_{1}>u_{2}\}}(X_{r})
(f_{1}(X_{r}, u_{1}(X_{r}))-f_{2}(X_{r}, u_{2}(X_{r})))\,dr \\&
\quad+ \int_{t}^{\tau}\mathbf{1}_{\{u_{1}>u_{2}\}}(X_{r})\,
d(A_{r}^{\mu_{1}}-A_{r}^{\mu_{2}})
+\int_{t}^{\tau}\mathbf{1}_{\{u_{1}>u_{2}\}}(X_{r})\,dA_{r}^{\nu_{1}}\\
&\quad-
\int_{t}^{\tau}\mathbf{1}_{\{u_{1}>u_{2}\}}(X_{r})\,dA_{r}^{\nu_{2}}
-\int_{t}^{\tau}\mathbf{1}_{\{u_{1}>u_{2}\}}(X_{r-})\,d(M_{r}^{1}-M_{r}^{2})\\
&=:\sum_{i=1}^{6}I_{i}(t,\tau).
\end{align*}
It is easy to see that $I_{2}(t,\tau)\le 0$ and $I_{3}(t, \tau)\le
0$. By the minimality of $\nu_1,\, \nu_2$ (condition (d) in the
definition of a solution of the obstacle problem), we have
\[
I_{4}(t, \tau)\le \int_{t}^{\tau}\mathbf{1}_{\{u_{1}>u_{2}\}}(u_{1}-u_{2})^{-1}
(u_{1}-u_{1}\wedge u_{2})\,dA_{r}^{\nu_{1}^{+}}=0
\]
and
\[
I_{5}(t, \tau)\le \int_{t}^{\tau}\mathbf{1}_{\{u_{1}>u_{2}\}}(u_{1}-u_{2})^{-1}
(u_{1}\vee u_{2}-u_{2})\,dA_{r}^{\nu_{2}^{-}}=0.
\]
The rest of the proof runs as in the proof of Proposition \ref{stw3.33}.
\end{proof}

\begin{corollary}
\label{wn4.44} Under \mbox{\rm{(H1)}} there exists at most one
solution of \mbox{\rm{OP}}$(f+d\mu, h_{1}, h_{2})$.
\end{corollary}
 Below we give the main theorem of this section.
We give an existence result for (\ref{eq1.1}) and show the convergence
of two penalization schemes. In the first one, we approximate the
solution $(u,\nu)$ to OP$(f+d\mu,h_1,h_1)$ by solutions $u_n$ to
PDE (\ref{eq4.00.1}) (with $n=k$). In (\ref{eq4.00.1}), a measure
$\nu_n$ with density (with respect to $m$) defined via $u_n$ and
$h_1,h_2$ appears. The convergence of $\nu_n$ to $\nu$ is in the
same metric as in the case of one barrier, i.e.
$A^{\nu_n}\rightarrow A^{\nu}$ in $\mathfrak S_c^q$. In the second
penalization scheme, we approximate $u$ by the first component of
the solution $(u_k,\alpha_k)$ to the obstacle problem
(\ref{eq4.00.2}) with one lower barrier $h_1$, and we approximate
$\nu$ by measures $\nu_k$ defined as the sum of  $\alpha_k$ and a
measure with density (with respect to $m$) defined via $ u_k,
h_2$. The advantage of the second penalization is that $\{u_k\}$
is monotone, and we have stronger convergence of the approximation
measures $\nu_k$ (see Corollary \ref{wn.wn}). As in the case of
one barrier, as a by-product we also get a probabilistic
representation of solutions.

\begin{theorem}
\label{th4.3}
Assume \mbox{\rm{(H1)--(H4)}}. Then there exists a
solution $(u,\nu)$ of \mbox{\rm{OP}}$(f+d\mu, h_{1},h_{2})$ if and only if
\mbox{\rm{(H6)}} is satisfied. Moreover, if \mbox{\rm(H6)} is
satisfied, then $u\in\mathfrak S^q_c$ for $q\in (0,1)$
and
\begin{enumerate}
\item[\rm(i)] if $u_{n,k}$ is a solution of the equation
\begin{equation}
\label{eq4.00.1}
-Lu_{n,k}=f(\cdot, u_{n,k})+\mu+n(u_{n,k}-h_{1})^{-}
-k(u_{n,k}-h_{2})^{+},
\end{equation}
then $u_{n,k}\rightarrow u$ q.e. and in
$\mathfrak S^q_c$ for $q\in (0,1)$ as $n,k\rightarrow\infty$,
\item[\rm(ii)]if $(u_{k}, \alpha_{k})$ is a solution of the
obstacle problem
\begin{equation}
\label{eq4.00.2}
-Lu_{k}=f(\cdot, u_{k})+\mu+\alpha_{k}
-k(u_{k}-h_{2})^{+}, \quad u_{k}\ge h_{1},
\end{equation}
then $u_{k}\searrow u$ q.e. and in $\mathfrak S^q_c$
for $q\in (0,1)$  as $k\rightarrow\infty$.
\end{enumerate}
\end{theorem}
\begin{proof}
The necessity is clear. To prove that
(H6) is sufficient let us first observe that  by Proposition
\ref{stw3.33}, $u_{k}\ge u_{k+1}$ and $d\alpha_{k}\le
d\alpha_{k+1}$. By (H6) there exist a function $v$ and a measure
$\lambda\in\mathbb M_0$ such that
\[
-Lv=\lambda, \quad f(\cdot, v)\in\mathbb M_0, \quad h_{1}\le v\le
h_{2} \quad m\rm{\mbox{-a.e.}}
\]
Hence
\[
-Lv=f(\cdot,v)+(\lambda^{+}+f^{-}(\cdot,v))
-(\lambda^{-}+f^{+}(\cdot,v))+
 n(v-h_{1})^{-}-k(v-h_{2})^{+}.
\]
Let $\overline{v}_{n}$ be a solution of the equation
\[
-L\overline{v}_{n}=f(\cdot,\overline{v}_{n})-\lambda^{-}-f^{+}(\cdot,v)-\mu^{-}
+n(\overline{v}_{n}-h_{1})^{-}.
\]
By Proposition \ref{stw3.33}, $\overline{v}_{n}\le v$ q.e., and consequently,
$\overline{v}_{n}\le h_{2}$, $m$-a.e.
Therefore
\[
-L\overline{v}_{n}=f(\cdot,\overline{v}_{n})-\lambda^{-}-f^{+}(\cdot,v)-\mu^{-}
+n(\overline{v}_{n}-h_{1})^{-}-k(\overline{v}_{n}-h_{2})^{+}.
\]
By Proposition \ref{stw3.33} again, $u_{n,k}\ge\overline{v}_{n}$
q.e., which implies that
\begin{equation}
\label{eq3.4}
n(u_{n,k}-h_{1})^{-}\le n(\overline{v}_{n}-h_{1})^{-}.
\end{equation}
By Theorem \ref{th.m1}, $\overline{v}_{n}\nearrow\overline{v}$ q.e.
where $(\overline{v}, \overline{\nu})$ is a solution of the
obstacle problem
\[
-L\overline{v}=f(\cdot,\overline{v})-\lambda^{-}-f^{+}(\cdot,v)-\mu^- +\overline{\nu},
\quad \overline{v}\ge h_{1}.
\]
Hence
\begin{equation}
\label{eq3.41} E_{x}\int_{0}^{\zeta}dA_{t}^{\overline{\nu}_{n}}
\rightarrow E_{x}\int_{0}^{\zeta}dA_{t}^{\overline{\nu}}
\end{equation}
for q.e. $x\in E$, where
$\overline{\nu}_{n}=n(\overline{v}_{n}-h_{1})^{-}\cdot m$. Write
$\alpha_{n,k}=n(u_{n,k}-h_{1})^{-}\cdot m$. By (\ref{eq3.4}),
$E_{x}\int_{0}^{\zeta}dA_{t}^{\alpha_{n,k}}\le
E_{x}\int_{0}^{\zeta}dA_{t}^{\overline{\nu}_{n}}$, whereas by
Theorem \ref{th.m1},
$E_{x}\int_{0}^{\zeta}dA_{t}^{\alpha_{n,k}}\rightarrow
E_{x}\int_{0}^{\zeta}dA_{t}^{\alpha_{k}}$ for q.e. $x\in E$, where $\alpha_k$ is defined in (ii).
Therefore
\begin{equation}
\label{eq3.5}
E_{x}\int_{0}^{\zeta}dA_{t}^{\alpha_{k}}\le
E_{x}\int_{0}^{\zeta}dA_{t}^{\overline{\nu}}
\end{equation}
for q.e. $x\in E$. Since $d\alpha_{k}\le d\alpha_{k+1}$,
\begin{equation}
\label{eq3.43}
dA_{t}^{\alpha_{k}}\le dA_{t}^{\alpha_{k+1}} \quad P_{x}\mbox{-a.s.}
\end{equation}
Set $A_{t}=\sup_{k\ge 1}A_{t}^{\alpha_{k}}$. By \cite[Lemma 3.2]{PengXu},
$A$ is a c\`adl\`ag process. Consequently, it is a positive additive functional  as an increasing   limit of additive functionals.
Thus, $w:=E_\cdot A_\zeta$ is an excessive function (see \cite[Proposition IV.2.4]{BG}). Consequently, by \cite[Theorem A.2.7]{Fukushima}, $w$
is finely-continuous. Therefore, by \cite[Theorem 4.6.1.]{Fukushima}, $w$ is quasi-continuous.
This implies that $A$ is a continuous AF.  Therefore there exists
a smooth measure $\alpha$ such that $A=A^{\alpha}$. Moreover, by
(\ref{eq3.41}) and (\ref{eq3.5}), $\alpha\in\mathbb M_0$. By
(\ref{eq3.43}) and Dini's theorem, for every $T>0$,
\begin{equation}
\label{eq3.52}
\sup_{t\le T}|A_{t}^{\alpha_{k}}-A_{t}^{\alpha}|\rightarrow_{P_{x}}0
\end{equation}
for q.e. $x\in E$. Let $u(x)=\inf_{k\ge 1}u_{k}(x)$, where $u_k$ is defined in (ii). Thanks
to (\ref{eq3.52}) we may now repeat arguments from   the
proof of Theorem \ref{th.m1} to show that $u$ is quasi-continuous,
and moreover, the following hold:
\[
E_{x}\int_{0}^{\zeta}|f_{u_{k}}(X_{t})-f_{u}(X_{t})|\,dt
\rightarrow 0
\]
for q.e. $x\in E$, there exists a nonnegative measure
$\delta\in\mathbb M_0$ such that for every $T>0$,
\begin{equation}
\label{eq3.6} \sup_{t\le
T}|A_{t}^{\delta_{k}}-A_{t}^{\delta}|\rightarrow_{P_{x}} 0
\end{equation}
for q.e. $x\in E$, where $\delta_{k}=k(u_{k}-h_{2})^{+}\cdot m$,
\begin{equation}
\label{eq3.7}
\sup_{t\le T}|u_{k}(X_{t})-u(X_{t})|\rightarrow_{P_{x}} 0
\end{equation}
for q.e. $x\in E$, and finally,
\begin{equation}
\label{eq4.122}
u(x)=E_{x}\int_0^{\zeta}f_u(X_t)\,dt+E_x\int_0^{\zeta}dA^{\mu}_t
+E_x\int_0^{\zeta}dA^{\alpha}_t-E_x\int_0^{\zeta}dA^{\delta}_t
\end{equation}
for q.e. $x\in E$. By (\ref{eq3.5}), $u\ge h_{1}$ $m$-a.e. By the
definition of a solution of the obstacle problem,
\[
u_k(x)=E_{x}\int_0^{\zeta}f_{u_k}(X_t)\,dt+E_x\int_0^{\zeta}dA^{\mu}_t
+E_x\int_0^{\zeta}dA^{\alpha_k}_t-E_x\int_0^{\zeta}dA^{\delta_k}_t
\]
for q.e. $x\in E$. From the above equation, (\ref{eq4.122}) and
the convergence results for  $u_{k}, f_{u_{k}}, A^{\alpha_{k}}$ we
have already proved,  we conclude that
\begin{equation}
\label{eq3.8} E_{x}\int_{0}^{\zeta}dA_{t}^{\delta_{k}}\rightarrow
E_{x}\int_{0}^{\zeta}dA_{t}^{\delta}
\end{equation}
for q.e. $x\in E$, which implies that $u\le h_{2}$ $m$-a.e. Using
(\ref{eq3.52})--(\ref{eq3.7}) we can show  in the same way as in
the proof of minimality of the measure $\nu$ in Theorem
\ref{th.m1} that for every quasi-continuous $h_{1}^{*}, h_{2}^{*}$
such that $h_{1}\le h_{1}^{*}\le u\le h_{2}^{*}\le h_{2}$ $m$-a.e.
we have
\[
E_{x}\int_{0}^{\zeta}(h_{2}^{*}(X_{t})-u(X_{t}))\,dA_{t}^{\delta}
=E_{x}\int_{0}^{\zeta}(u(X_{t})-h_{1}^{*}(X_{t}))\,dA_{t}^{\alpha}=0
\]
for q.e. $x\in E$. Of course, putting $\nu=\delta-\alpha$ yields
the  above equation with $\nu^{-}$ in place of $\delta$ and
$\nu^{+}$ in place of $\alpha$. Thus the pair $(u,\nu)$ is a
solution of OP$(f+d\mu,h_1,h_2)$. Observe that
\begin{equation}
\label{eq3.9} w_{n}\le u_{n,k}\le u_{k} \quad \mbox{q.e.},
\end{equation}
where $(w_{n}, \beta_{n})$ is a solution of the obstacle problem
\[
-Lw_{n}=f(\cdot, w_{n})+ n(w_{n}-h_{1})^{-}+\mu-\beta_{n},
\quad w_{n}\le h_{2}.
\]
To see this it is enough to observe that
\[
-Lu_{k}=f(\cdot, u_{k})+ n(u_{k}-h_{1})^{-}-k(u_{k}-h_{2})^{+}+ \mu+\alpha_{k}
\]
and
\[
-Lw_{n}=f(\cdot, w_{n})+ n(w_{n}-h_{1})^{-}-k(w_{n}-h_{2})^{+}+ \mu-\beta_{n},
\]
and apply Proposition \ref{stw3.33}. By the same method as in the
case of $\{u_k\}$, one can show that the limit of $\{w_n\}$ is the
first component of the solution of OP$(f+d\mu,h_1,h_2)$. Hence, by
Corollary \ref{wn4.44}, $w_{n}\rightarrow u$ q.e. Finally, observe
that by (\ref{eq3.5}) and (\ref{eq3.52})--(\ref{eq3.7}), for every
$q\in (0,1)$,
\[
E_{x}\sup_{t\ge 0}|A_{t}^{\delta_{k}}-A_{t}^{\delta}|^{q}
+E_{x}\sup_{t\ge 0}|A_{t}^{\alpha_{k}}-A_{t}^{\alpha}|^{q}
\rightarrow 0
\]
for q.e. $x\in E$. Moreover, by the Tanaka-Meyer formula (see
\cite[Theorem IV.70]{Protter}),
\[
|u_{k}(X_{t})|\le E_{x}\Big(\int_{0}^{\zeta}|f(X_{t}, 0)|\,dt
+\int_{0}^{\zeta} dA_{t}^{|\mu|}+ \int_{0}^{\zeta}dA_{t}^{|\nu|}|
\mathcal{F}_{t}\Big).
\]
Therefore by \cite[Lemma 6.1]{BDHPS}, for every $q\in (0,1)$,
\[
E_{x}\sup_{t\ge 0}|u_{k}(X_{t})|^{q} \le
(1-q)^{-1}\Big[E_{x}\Big(\int_{0}^{\zeta}|f(X_{t},0)|
+\int_{0}^{\zeta}dA_{t}^{|\mu|}+
\int_{0}^{\zeta}dA_{t}^{|\nu|}\Big)\Big]^{q}.
\]
From this we conclude that $u_{n}\rightarrow u$ in
$\mathfrak S^q_c$ for $q\in (0,1)$. In the same manner
we can see that $w_{n}\rightarrow u$ in
$\mathfrak S^q_c$ for $q\in(0,1)$, which when combined
with (\ref{eq3.9}) implies that $u_{n,k}\rightarrow u$ in
$\mathfrak S^q_c$ for $q\in(0,1)$.
\end{proof}

\begin{corollary}
\label{wn.wn} Assume \mbox{\rm(H1)--(H4), (H6)} and retain the
notation from Theorem \ref{th4.3} and its proof. Then for every
$q\in(0,1)$,  and for q.e. $x\in E$,
\begin{enumerate}
\item[\rm(i)]
$ E_{x}\sup_{t\ge 0}|A_{t}^{\alpha_{k}}-A_{t}^{\nu^{+}}|^{q}+
E_{x}\sup_{t\ge 0}|A_{t}^{\delta_{k}}-A_{t}^{\nu^{-}}|^{q}
\rightarrow 0$ as $k\rightarrow \infty $,
\item[\rm(ii)] $E_{x}\sup_{t\ge 0}|A_{t}^{\nu_{n}
}-A_{t}^{\nu}|^{q} \rightarrow 0$ as $n\rightarrow \infty$, where
$\nu_n=n(u_{n,n}-h_1)^--n(u_{n,n}-h_2)^+$.
\end{enumerate}
\end{corollary}
\begin{proof}
(i) One can regard $(u,\nu^{-})$ as a
solution of OP$(f+d\mu+d\nu^{+}, h_{2})$ (with upper barrier).
Therefore by Theorem \ref{th.m1}, $y_{k}\searrow u$ q.e., where
\[
-Ly_{k}=f(\cdot, y_{k})-k(y_{k}-h_{2})^{+}+\nu^{+}+\mu,
\]
and for every $q\in(0,1)$,
\begin{equation}
\label{eq3.10}
E_{x}\sup_{t\ge 0}|A_{t}^{\beta_{k}}-A_{t}^{\nu^{-}}|^{q}
\rightarrow 0
\end{equation}
for q.e. $x\in E$, where $\beta_{k}=k(y_{k}-h_{2})^{+}\cdot m$.
Since $y_{k}\searrow u$,  $y_{k}\ge h_{1}$ q.e. Therefore
\[
-Ly_{k}=f(x, y_{k})+n(y_{k}-h_{1})^{-}
-k(y_{k}-h_{2})^{+}+\nu^{+}+\mu.
\]
By Proposition \ref{stw3.33}, $y_{k}\ge u_{n,k}$ q.e., so
$k(u_{n,k}-h_{2})^{+}\le k(y_{k}-h_{2})^{+}$. By (\ref{eq3.10})
and the convergence of $\{A^{\alpha_{n,k}}\}$ showed in the proof
of Theorem \ref{th.m1}, $dA^{\alpha}\le dA^{\nu^{+}}$,  which
implies that $d\alpha\le d\nu^{+}$. The same reasoning applied to
the measure $\delta$ shows that $d\delta \le d\nu^{-}$. From this
and minimality of the Jordan  decomposition of measure $\nu$ we
conclude that $\alpha=\nu^{+}$, $\delta=\nu^{-}$.

 (ii) By Theorem \ref{th4.3}, $u_{n,n}\rightarrow
u$ in $\mathfrak S^q_c$ for every $q\in (0,1)$. By (\ref{eq3.9}),
$w_1\le u_{n,n}\le u_1,\, n\ge 1$. The rest of the proof of (ii)
is analogous to that of Corollary \ref{wn.wn.wn.wn}.
\end{proof}

\begin{proposition}
\label{stw.r2} Assume that \mbox{\rm{(H1)--(H4), (H6)}} are satisfied and let
$(u,\nu)$ be a solution of \mbox{\rm{OP}}$(f+d\mu,h_1,h_2)$. Then
\[
\|\nu^+\|\le 4(\|\mu\|+\|f_0\|+\|\lambda^+\|+\|f^-_v\|)
\]
and
\[
\|\nu^-\|\le 4(\|\mu\|+\|f_0\|+\|\lambda^-\|+\|f^+_v\|)
\]
with $\lambda=-Lv$, where $v$ is the function from condition
\mbox{\rm(H6)}.
\end{proposition}
\begin{proof}
From (\ref{eq3.5}), (\ref{eq3.43}) and \cite[Lemma 2.6]{KR:CM} we
deduce that $\|\alpha\|\le \|\overline{\nu}\|$. Hence
$\|\nu^+\|\le \|\overline{\nu}\|$ since $\alpha=\nu^+$ by
Corollary \ref{wn.wn}. On the other hand, by Proposition
\ref{stw.r1},
\[
\|\overline{\nu}\|\le 2(\|\lambda^+\|+\|f^-_{v}\|+\|\mu^-\|
+\|f_0\|+\|\lambda^+\|+\|f_v^-\|),
\]
which proves the desired inequality for $\nu^+$. The inequality
for $\nu^-$ can be proved in much the same way.
\end{proof}

\begin{proposition}
Assume that \mbox{\rm{(H1)--(H4), (H6)}}are satisfied  and let $(u,\nu)$ be a
solution of \mbox{\rm{OP}}$(f+d\mu,h_1,h_2)$. If
$\lambda,f_v,f_0,\mu\in\mathcal{M}_{0,b}$, then
$\nu\in\mathcal{M}_{0,b}$, $T_{k}(u)\in D_e[\EE]$ for every $k\ge0$
and \mbox{\rm(\ref{eq3.13})} is satisfied.
\end{proposition}
\begin{proof}
Follows from Proposition \ref{stw.r2} and \cite[Proposition 3.10,
Theorem 4.2]{KR:CM}.
\end{proof}

\begin{proposition}
\label{stw.stw123}
Let $(u,\nu)$ be a solution to \mbox{\rm OP}$(d\mu,h_1,h_2)$. Assume that there exists $v$ such that
$h_1\le v\le h_2$ and $v=R\lambda$ for some $\lambda$ such that $|\lambda|\in D'_e[\EE]$. Then
$u\in D_e[\EE]$, $\nu\in D'_e[\EE]$ and
$(u,\nu)$ is the unique pair in $D_e[\EE]\times D_e'[\EE]$ such that
\begin{equation}
\label{eq3.varne}
\EE(u,\eta)=\int_E\eta\,d\mu+\int_E\eta\,d\nu,\quad \eta\in D_e[\EE],\qquad h_1\le u\le h_2\,\,\mbox{a.e.}
\end{equation}
and
\[
 \int_E(u-\eta)\,d\nu\le 0,\quad \eta\in D_e[\EE],\,h_1\le \eta \le h_2\,\,\mbox{a.e.}
\]
\end{proposition}
\begin{proof}
Since $|\lambda|\in D'_e[\EE]$, $v\in D_e[\EE]$. With the notation of Theorem \ref{th4.3} (with $f\equiv 0$),
we have by (\ref{eq3.4}) that
\[
\|\alpha_{n,k}\|_{\EE'}\le\|\bar\nu_n\|_{\EE'},\quad n,k\ge1.
\]
Of course $(\bar u_n,\bar \nu_n)$ is a solution to OP$(-d\mu^--d\lambda^-,h-(u_n-h)^-)$, so by Proposition \ref{stw3.var},
\[
\|\bar\nu_n\|_{\EE'}\le3(\|\mu^-\|_{\EE'}+\|\lambda^-\|_{\EE'}+\|v\|_{\EE}).
\]
Since $\|R\beta\|_\EE\le\|\beta\|_{\EE'}$ for every $\beta\in D'_e[\EE]$, from the above inequalities it follows that
\[
\|R\alpha_{n,k}\|_{\EE}\le 3(\|\mu^-\|_{\EE'}+\|\lambda^-\|_{\EE'}+\|v\|_{\EE}).
\]
By Theorem \ref{th4.3} and Corollary \ref{wn.wn}, $R\alpha_{n,k}\rightarrow R\alpha_n$ as $k\rightarrow\infty$ and $R\alpha_n\nearrow R\nu^+$ as $n\rightarrow\infty$.
Hence we get
\[
\|R\nu^+\|_{\EE}\le 3(\|\mu^-\|_{\EE'}+\|\lambda^-\|_{\EE'}+\|v\|_{\EE}).
\]
This implies that $\nu^+\in D'_e[\EE]$.
Of course $(-u,\nu^-)$ is a solution to OP$(-d\mu-d\nu^+,-h_2)$, so
the desired result follows from
Proposition \ref{stw3.var}.
\end{proof}

\begin{proposition}
Assume \mbox{\rm{(H1)--(H4)}}. If $(u,\nu)$ is a solution of
\mbox{\rm{OP}}$(f+d\mu, h_{1}, h_{2})$, then $u$ admits
representation \mbox{\rm(\ref{eq1.4})}.
\end{proposition}
\begin{proof}
Let $v$ be a supersolution of PDE$(f+d\mu-d\nu^{-})$ such that
$v\ge h_{1}$ $m$-a.e. Then there exists a nonnegative measure
$\lambda\in\mathbb M_0$ such that
\[
-Lv=f(x,v)+\mu-\nu^{-}+\lambda.
\]
Since $v\ge h_{1}$ $m$-a.e.,
\[
-Lv=f(x,v)+\mu+ n(v-h_{1})^{-}-\nu^{-}+\lambda.
\]
Observe that the pair $(u, \nu^{+})$ is a solution of
OP$(f+d\mu-d\nu^{-}, h_{1})$. Therefore, by Theorem \ref{th.m1},
$u_{n}\nearrow u$ q.e., where
\[
-Lu_{n}=f(x,u_{n})+\mu+ n(u_{n}-h_{1})^{-}-\nu^{-}.
\]
By Proposition \ref{stw3.33}, $u_{n}\le v$ q.e., which implies
that $u\le v$ q.e.
\end{proof}

\section{Lewy-Stampacchia type inequality and stability results}
\label{sec5}

In this section, we prove  Lewy-Stampacchia type inequality in our
general framework and give some stability results for solutions.
In the case of one barrier and regular data, inequalities of such
type for nonlocal operators (on $\mathbb{R}^{n}$) were proved  in
\cite{SV} (see also the recent papers
\cite{GM,PV} for abstract Lewy-Stampacchia inequality and for the
same type of inequality in the Heisenberg group).

Let us stress that the measures $f_{h_1}\cdot m,\mu,Lh_1$ and
$\nu$ in the theorem below need not be finite.

\begin{theorem}
Let $\mu\in\mathbb M_0$ and let $(u,\nu)$ be a solution of
\mbox{\rm{OP}}$(f+d\mu, h_{1}, h_{2})$. If $h_{1}$ is a difference
of natural potentials, then
\begin{equation}
\label{eq4.10}
\nu^{+}\le \mathbf{1}_{\{u=h_{1}\}}\cdot (f_{h_{1}}{\cdot m}+\mu+Lh_{1})^{-}.
\end{equation}
\end{theorem}
\begin{proof}
By the assumption on the barrier $h_1$,  there exists a measure $\alpha\in\mathbb M_0$ such that for q.e. $x\in E$,
\[
h_{1}(x)=E_{x}\int_{0}^{\zeta}dA_{t}^{\alpha}.
\]
Therefore, by Lemma \ref{lm3.1}, there exists a martingale AF $M^{1}$ of
$\mathbf{M}$ such that
\[
h_{1}(X_{t})=\int_{t}^{\zeta}dA_{t}^{\alpha}-\int_{t}^{\zeta}dM_{t}^{1},
\quad t\in[0,\zeta].
\]
By the Tanaka-Meyer formula (see \cite[Theorem IV.70]{Protter}),
\begin{align*}
(u-h_{1})^{+}(X_{t})&=(u-h_{1})^{+}(X_{0})-
\int_{0}^{t}\mathbf{1}_{\{u>h_{1}\}}(X_{r})f_{u}(X_{r})\,dr \\
&\,\,-\int_{0}^{t}\mathbf{1}_{\{u>h_{1}\}}(X_{r})\,
d(A_{r}^{\nu^{+}}+A_{r}^{\mu}-A_{r}^{\alpha})+
\int_{0}^{t}\mathbf{1}_{\{u>h_{1}\}}(X_{r})\,dA_{r}^{\gamma^{-}} \\
&\,\,-\frac 1 2 L_{t}^{0}(Y)+J_{t}^{+}+
\int_{0}^{t}\mathbf{1}_{\{u>h_{1}\}}(X_{r})\,d(M_{r}-M_{r}^{1}),
\end{align*}
where
\[
J_{t}^{+}=\sum_{0<s\le t}(\varphi(Y_{s})-\varphi(Y_{s-})
-\varphi'(Y_{s-})\Delta Y_{s}),\quad
Y_{t}=(u-h_{1})(X_{t}),\quad \varphi(x)=x^{+},
\]
$\varphi'$ denotes the left derivative of $\varphi$, and $ L^{0}(Y)$ is the local time of $Y$ at $0$. Since
$Y_{t}\ge 0$, $t\ge 0$, we conclude from the above equations that
\begin{align*}
0&= \int_{0}^{t}\mathbf{1}_{\{u=h_{1}\}}(X_{r})f_{h_{1}}(X_{r})\,dr +
\int_{0}^{t}\mathbf{1}_{\{u=h_{1}\}}(X_{r})\,
d(A_{r}^{\nu^{+}}+A_{r}^{\mu}-A_{r}^{\alpha}) \\
&\,\,
-\int_{0}^{t}\mathbf{1}_{\{u=h_{1}\}}(X_{r})\,dA_{r}^{\nu^{-}}
+\frac 1 2 L_{t}^{0}(Y)+J_{t}^{+}
-\int_{0}^{t}\mathbf{1}_{\{u=h_{1}\}}(X_{r-})\,d(M_{r}-M_{r}^{1}).
\end{align*}
Since $\int_{0}^{t}dA_{r}^{\nu^{+}}
=\int_{0}^{t}\mathbf{1}_{\{u=h_{1}\}}(X_{r})\,dA_{r}^{\nu^{+}}$,
\begin{align*}
\frac 1 2 L_{t}^{0}(Y)+J_{t}^{+,p}+\int_{0}^{t}dA_{r}^{\nu^{+}} &=
-\int_{0}^{t}\mathbf{1}_{\{u=h_{1}\}}(X_{r})f_{h_{1}}(X_{r})\,dr
\\&
\quad+\int_{0}^{t}\mathbf{1}_{\{u=h_{1}\}}(X_{r})\,
d(A_{r}^{\nu^{-}}-A_{r}^{\mu}+A_{r}^{\alpha}),
\end{align*}
where $J_{t}^{+,p}$ is the dual predictable projection of the
process $J_{t}^{+}$. Since $dA^{\nu^{+}}, dA^{\nu^{-}}$ are
orthogonal, $\int_{0}^{t}\mathbf{1}_{\{u=h_{1}\}}(X_{r})\,dA_{r}^{\nu^{-}}=0$. Therefore
\begin{align*}
dA_{t}^{\nu^{+}}& \le \mathbf{1}_{\{u=h_{1}\}}(X_{t})(-f_{h_{1}}(X_{t})\,dt-
dA_{t}^{\mu}+dA_{t}^{\alpha})^{+} \\
&=\mathbf{1}_{\{u=h_{1}\}}(X_{t})(f_{h_{1}}(X_{t})\,dt+
dA_{t}^{\mu}-dA_{t}^{\alpha})^{-},
\end{align*}
which combined with Revuz duality implies (\ref{eq4.10}).
\end{proof}

\begin{proposition}
\label{stw.stab}
Assume that $\mu_n,\mu\in\mathbb M_0$ and  $f_n, f$ satisfy \mbox{\rm(H1)}. Let $(u_n,\nu_n)$, $(u,\nu)$ be solutions of \mbox{\rm OP}$(f_n+d\mu_n,h_1,h_2)$ and
\mbox{\rm OP}$(f+d\mu,h_1,h_2)$, respectively. If
\begin{equation}
\label{eq5.2}
R|\mu_n-\mu|\rightarrow 0,\qquad R|f_n(\cdot,u)-f(\cdot,u)|\rightarrow 0\quad m\mbox{-a.e.},
\end{equation}
then $u_n\rightarrow u$ $m$-a.e.
\end{proposition}
\begin{proof}
By the definition of a solution to the obstacle problem and  Lemma \ref{lm3.1}, there exist  martingale AFs $M,M^n$ of $\mathbf{M}$
such that for q.e. $x\in E$,
\[
u(X_{t})=\int_{t}^{\zeta}f(X_r,u(X_r))\,dr
+\int_{t}^{\zeta}dA_{r}^{\mu}+\int_{t}^{\zeta}dA_{r}^{\nu} +\int_{t}^{\zeta}dM_{r},
\quad 0\le t\le\zeta,
\]
and
\[
u_n(X_{t})=\int_{t}^{\zeta}f_n(X_r,u_n(X_r))\,dr
+\int_{t}^{\zeta}dA_{r}^{\mu_n}+\int_{t}^{\zeta}dA_{r}^{\nu_n} +\int_{t}^{\zeta}dM^n_{r},
\,\, 0\le t\le\zeta,
\]
$P_x$-a.s. By the Tanaka-Meyer formula, (H1) and the minimality conditions for $\nu_n$ and $\nu$ we have
\begin{align*}
|u_n(x)-u(x)|&\le E_x\int_0^\zeta|f_n(\cdot,u)-f(\cdot,u)|(X_r)\,dr+E_x\int_0^\zeta dA^{|\mu_n-\mu|}_r\\
&=R|f_n(\cdot,u)-f(\cdot,u)|(x)+R|\mu_n-\mu|(x)
\end{align*}
for q.e. $x\in X$.  By this and (\ref{eq5.2}), $u_n\rightarrow u$ $m$-a.e.
\end{proof}

\begin{remark}
\label{rem.uts}
If $\mu_n\rightarrow \mu$ in the total variation norm and $f_n(\cdot,u)\rightarrow f(\cdot,u)$ in $L^1(E;m)$, then assumption (\ref{eq5.2}) is satisfied for some subsequence of $\{n\}$.
Indeed, since $\EE$ is
transient, there exists a strictly positive $\eta\in\BB_b(E)$ such that $\|\hat G\eta\|_\infty<\infty$ (see \cite[Theorem 1.3.4]{O}). Therefore
\[
\int_E\eta R|\mu_n-\mu|\le \|\hat G\eta\|_\infty|\mu_n-\mu|(E),
\]
and
\[
\int_E \eta R|f_n(\cdot,u)-f(\cdot,u)|\le \|\hat G\eta\|_\infty \|f_n-f\|_{L^1},
\]
from which the desired result follows.
\end{remark}

\begin{remark}
Let $(u,\nu)$ be a solution to \mbox{\rm OP}$(f+d\mu,h_1,h_2)$. Assume that there exists $v$ such that
$h_1\le v\le h_2$ and $v=R\lambda$ for some $\lambda$ such that $|\lambda|\in D'_e[\EE]$ (in the case where $h_2\equiv +\infty$ it is enough to
assume that there exists $v\in D_e[\EE]$ such that $v\ge h_1$). Let $g$ be a strictly positive function such that $g\in D_e'[\EE]$ and let $\{F_n\}$ be a nest such that $\mu_n:=\mathbf 1_{F_n}\cdot \mu\in D'_e[\EE]$. For $n\in\BN$ set
\[
f_n(x,y)=\frac{ng(x)}{1+ng(x)}(f\wedge n)(x,y), \quad x\in E,\,y\in\BR.
\]
By Theorem \ref{th4.3} and Proposition \ref{stw.stw123}, there exists a unique solution $u_n$ of variational inequality (\ref{eq1.1234}) with $f,\mu$ replaced by $f_n,\mu_n$, and moreover, $u_n$ coincides with solution to OP$(f_n+d\mu_n,h_1,h_2)$. By Proposition \ref{stw.stab} and Remark \ref{rem.uts}, up to subsequence, $u_n\rightarrow u$ $m$-a.e. This shows that  each solution to (\ref{eq1.1}) may be approximated by solutions to variational inequalities.
\end{remark}

\subsection*{Acknowledgements}
The research was supported by Polish National Science Centre (Grant No.
2012/07/D/ST1/02107).


\end{document}